\documentclass[11pt]{article}

\usepackage{amsmath}
\usepackage{amssymb}
\usepackage{theorem}
\usepackage{euscript}
\usepackage{amscd}

\renewcommand{\title}[1]{
     \addvspace{3\baselineskip}  
     \begin{center} \LARGE \bf #1
     \end{center}
     \addvspace{2\baselineskip}}   

\renewcommand{\author}[1]{
     \addvspace{-1\baselineskip}  
     \begin{center} \large \sc #1
     \end{center}
     \addvspace{2\baselineskip}}   

\makeatletter

\def\section{%
        \@startsection{section}{1}{\z@}%
        {8ex plus 6ex minus 3ex}{\baselineskip}%
        {\normalfont\large\scshape\centering}%
        }

\renewcommand{\paragraph}[1]{{\par\removelastskip\vskip.5\baselineskip%
         \indent{\itshape{#1}}{\ifperiod.\else\global\periodtrue\fi}%
         \rm \ignorespaces}}

\let\goth=\mathfrak
\let\calligraphy=\mathcal

\def\CC{{\mathbb C}}

\def\NN{{\mathbb N}}

\def\PP{{\mathbb P}}
\def\QQ{{\mathbb Q}}
\def\RR{{\mathbb R}}
\def\SS{{\mathbb S}}

\def\ZZ{{\mathbb Z}}

\def\Aa{{\calligraphy A}}
\def\Bb{{\calligraphy B}}

\def\Ii{{\calligraphy I}}
\def\Jj{{\calligraphy J}}

\def\Ll{{\calligraphy L}}

\def\Oo{{\calligraphy O}}
\def\Pp{{\calligraphy P}}

\def\Uu{{\calligraphy U}}

\def\aaa{{\goth a}}
\def\BBB{{\goth B}}
\def\bbb{{\goth b}}

\def\FFF{{\goth F}}

\def\mmm{{\goth m}}

\def\qqq{{\goth q}}
\def\RRR{{\goth R}}

\def\VVV{{\goth V}}

\def\hbar{{\,\overline{\!h}}}

\def\Mbar{{\,\overline{\!M}}}

\def\Ctilde{{\,\widetilde{\!C}}}

\def\Htilde{{\,\widetilde{\!H}}}

\def\Ghat{{\,\widehat{\!G}}}

\def\Hhat{{\,\widehat{\!H}}}

\def\cC{{\boldsymbol{C}}}
\def\dd{{\boldsymbol{d}}}

\def\ee{{\boldsymbol{e}}}
\def\eE{{\boldsymbol{E}}}

\def\ll{{\boldsymbol{l}}}

\def\nn{{\boldsymbol{n}}}

\def\xx{{\boldsymbol{x}}}
\def\yy{{\boldsymbol{y}}}

\def\zero{{\boldsymbol{0}}}

\def\card{\operatorname{card}}

\def\except{\operatorname{except}}

\def\Hom{\operatorname{Hom}}

\def\Im{\operatorname{Im}}

\def\Res{\operatorname{Res}}

\def\Specbf{\mathbf{Spec\,}}

\let\lra=\longrightarrow

\def\ie{{\it i.e.}~}

\def\inv{^{-1}}

\let\phi=\varphi
\let\epsilon=\varepsilon

\newcommand{\floor}[1]{\left\lfloor#1\right\rfloor}
\newcommand{\ceil}[1]{\left\lceil#1\right\rceil}
\newcommand{\fpart}[1]{\left\langle#1\right\rangle}

\newcommand{\textq}[1]{\quad\text{#1}\quad}
\newcommand{\textql}[1]{\quad\text{#1}}
\newcommand{\textqr}[1]{\text{#1}\quad}

\newcommand{\listspace}{\setlength{\itemsep}{-.5pt}}

\def\theoname{Theorem}
\def\lemmaname{Lemma}
\def\propositionname{Proposition}
\def\notationname{Notation}
\def\corollaryname{Corollary}
\def\conjecturename{Conjecture}
\def\remarkname{Remark}
\def\remarksname{Remarks}
\def\examplename{Example}
\def\examplesname{Examples}
\def\definitionname{Definition}
\def\definitionsname{Definitions}
\def\notationname{Notation}

\def\proofname{Proof}

\def\Dquad{\hskip 0.6em plus .02em minus .2em}  
\def\Dpar{\belowdisplayskip=0pt\belowdisplayshortskip=0pt\par}

\def\bigpenalty{\interlinepenalty=\@M}
\def\smallpenalty{\interlinepenalty=100}

\newif\ifperiod \periodtrue 

\def\D@makemargins{%
  \labelsep=0pt
  \itemindent=0pt
  \labelwidth=0pt}

\def\D@restoremargins{%
  \labelsep=5pt
  \itemindent=0pt
  \leftmargin=5mm  
  \labelwidth=\leftmargin \advance\labelwidth by -\labelsep}

\def\th@Dindent{\hspace\parindent}
\def\th@Dheadingshape{\scshape}

\gdef\th@DthAndSuchtheo{%
  \D@makemargins%
  \def\@begintheorem##1##2{%
  \item[]\th@Dindent{\th@Dheadingshape ##1~\rm ##2.}\Dquad         
        \D@restoremargins}%
  \def\@opargbegintheorem##1##2##3{\def\next{##3}%
  \item[]\th@Dindent{\th@Dheadingshape ##1~\rm ##2\ifx\next\empty
  \else\ {\normalfont(##3)}\fi.}         
        \D@restoremargins}}

\gdef\th@DthAndSuchtheostar{%
  \D@makemargins%
  \def\@begintheorem##1##2{%
  \item[]\th@Dindent{\th@Dheadingshape ##1.}\Dquad     
        \D@restoremargins}%
  \def\@opargbegintheorem##1##2##3{\def\next{##3}%
  \item[]\th@Dindent{\th@Dheadingshape ##1\ifx\next\empty
  \else\ ##3\fi.}\Dquad         
        \D@restoremargins}}

\gdef\th@DthAndSuchliketheo{
  \D@makemargins%
  \def\@begintheorem##1##2{%
    \@latex@error{likethm: You must provide an argument in square brackets,
    though it may be empty [] !}%
    }%
  \def\@opargbegintheorem##1##2##3{%
        \def\next{##3}\ifx\next\empty\item[\th@Dindent]\else
        \item[]\th@Dindent{\th@Dheadingshape \next.}\Dquad\fi
        \D@restoremargins}}

\gdef\th@DdefAndSuch{%
  \D@makemargins%
  \def\@begintheorem##1##2{%
  \item[]\th@Dindent{\def@Dheadingshape ##1~\rm ##2.}\Dquad         
        \D@restoremargins}%
  \def\@opargbegintheorem##1##2##3{\def\next{##3}%
  \item[]\th@Dindent{\def@Dheadingshape ##1~\rm ##2\ifx\next\empty
  \else\ {\normalfont(##3)}\fi.}         
        \D@restoremargins}}

\gdef\th@DdefAndSuchStar{%
  \D@makemargins%
  \def\@begintheorem##1##2{%
  \item[]\th@Dindent{\def@Dheadingshape ##1.}\Dquad     
        \D@restoremargins}%
  \def\@opargbegintheorem##1##2##3{\def\next{##3}%
  \item[]\th@Dindent{\th@Dheadingshape ##1\ifx\next\empty
  \else\ ##3\fi.}\Dquad         
        \D@restoremargins}}

\def\th@Dheadingshape{\scshape}
\def\def@Dheadingshape{\itshape}

\theoremstyle{DthAndSuchliketheo}
\theorembodyfont{\bigpenalty\itshape}   
\newtheorem{likethm}{}
\theorembodyfont{\rmfamily}  

\theoremstyle{DthAndSuchtheostar}
\theorembodyfont{\bigpenalty\itshape}
\newtheorem{thm*}{\theoname}
\newtheorem{lem*}{\lemmaname}
\newtheorem{pro*}{\propositionname}
\newtheorem{cor*}{\corollaryname}
\newtheorem{conjecture*}{\conjecturename}
\theorembodyfont{\smallpenalty\rmfamily}
\newtheorem{notation*}{\notationname}
\newtheorem{exa*}{\examplename}
\newtheorem{examples*}{\examplesname}
\theoremstyle{DdefAndSuchStar}
\theorembodyfont{\smallpenalty\rmfamily}
\newtheorem{definition*}{\definitionname}
\newtheorem{definitions*}{\definitionsname}
\newtheorem{rem*}{\remarkname}
\newtheorem{remarks*}{\remarksname}
\theoremstyle{DthAndSuchtheo}
\theorembodyfont{\bigpenalty\itshape}
\newtheorem{thm}{\theoname}[section]
\newtheorem{lem}[thm]{\lemmaname}
\newtheorem{pro}[thm]{\propositionname}
\newtheorem{cor}[thm]{\corollaryname}

\theorembodyfont{\smallpenalty\rmfamily}
\newtheorem{exa}[thm]{\examplename}

\newtheorem{notation}[thm]{\notationname}
\theoremstyle{DdefAndSuch}
\theorembodyfont{\smallpenalty\rmfamily}
\newtheorem{definition}[thm]{\definitionname}
\newtheorem{rem}[thm]{\remarkname}
\newtheorem{definitions}[thm]{\definitionsname}
\newtheorem{remarks}[thm]{\remarksname}

\theoremstyle{DthAndSuchtheo}               
\theorembodyfont{\itshape}

\newcommand{\proof}[1][]{{\par\removelastskip\vskip.6\baselineskip   
    \noindent\th@Dindent\def\next{#1}%
    {\itshape\proofname\ifx\next\empty\else\next\fi\ifperiod.%
      \else\global\periodtrue\fi\Dquad}%
    \clubpenalty=5000\rm\ignorespaces}\setcounter{step}{0}}

\newcounter{step}

\newcommand{\likeproof}[1][]{{\par\removelastskip\vskip.6\baselineskip
    \noindent\th@Dindent\def\next{#1}%
    {\itshape\ifx\next\empty\else\next\fi\ifperiod.%
      \else\global\periodtrue\fi\Dquad}%
    \clubpenalty=5000\rm\ignorespaces}\hspace{-2pt}\setcounter{step}{0}}

\def\qed{{\ifmmode\hskip 6mm plus 1mm minus 3mm{$\square$}
    \else
    \hfil\penalty50\hskip1em\null\nobreak\hfil
    {\hfill $\square$\parfillskip=0pt\finalhyphendemerits=0
      \let\par=\endgraf\par}
    \fi
    \Dpar\penalty-150\vskip.6\normalbaselineskip}}

\makeatother

\usepackage{pstricks}
\usepackage{pst-node,pst-tree}
\usepackage[arrow,curve,matrix,tips,frame]{xy}
\newpsobject{showgrid}{psgrid}{subgriddiv=1,griddots=10,gridlabels=6pt}

\begin{document}

\title{Mixed multiplier ideals and the irregularity of abelian
  coverings of the projective plane} 
\author{Daniel Naie}


\bigskip

\begin{abstract}
\noindent
A formula for the irregularity of abelian coverings of the projective
plane is established and some applications are presented.
\end{abstract}

\section{Introduction}

The initial intent of this study was to extend the formula for the
cyclic multiple planes from \cite{Na} to the case where the branching
curve $C$ is not transverse to the line at infinity $H_\infty$.  In
the transverse case, if $S$ denotes a desingularization of the
$\ZZ/n\ZZ$-cyclic covering of the plane associated to $C$ and
$H_\infty$, then
\begin{equation}
  \label{eq:q_cC}
  q(S) = \sum_{\substack{\xi\text{ jumping number of }C\\
      \xi\in1/(n\wedge\deg C)\,\ZZ,\,\,0<\xi<1}}
  h^1(\PP^2,\Oo_{\PP^2}(-3+\xi\deg C)\otimes\Jj(\xi\cdot C)).
\end{equation}
Hence, the irregularity is quasi-constant as a function of $n$, unlike
what happens in the non transverse case, when, as we see in
Example~\ref{ex:cyclic}, the irregularity might be a degree $1$
quasi-polynomial of $n$.  To understand the difference and to extend
the above formula to the non transverse case, we consider abelian
instead of cyclic coverings.  The role played by the multiplier ideals
will be taken by the mixed multiplier ideals.  Consequently, the goal
of this paper is to apply the theory of mixed multiplier ideals to
compute the irregularity of the abelian coverings of the projective
plane.

If $X$ is a smooth surface and $\aaa_1,\ldots,\aaa_t\subset\Oo_X$ are
non-zero ideal sheaves, the mixed multiplier ideal
$\Jj(\aaa_1^{x^1}\cdots\aaa_t^{x^t})$ varies with the rational vector
$\xx=(x^1,\ldots,x^t)\in\RR_+^t$.  Proposition \ref{pr:mainTool} and
Proposition \ref{pr:mainR} assert that there is a set of hyperplanes
called {\it jumping walls} with the following properties:
\begin{enumerate}\listspace
\item 
If the mixed multiplier ideal jumps, then the vector $\xx$ crosses a
jumping wall.  Consequently, the fibres of the map
$\xx\mapsto\Jj(\aaa_1^{x^1}\cdots\aaa_t^{x^t})$ are finite unions of
rational convex polytopes cut out by the jumping walls.
\item
The jumping walls are determined by the {\it jumping numbers} of the {\it
simple complete relevant ideals} (see Definition~\ref{d:relIdeals})
associated to the ideals $\aaa_i$.
\end{enumerate}
These results together with O.~Zariski's original idea introduced in
\cite{Za2} enable us to generalize formula (\ref{eq:q_cC}) to abelian
coverings of the projective plane.  Such a covering induces a
partition of the branching curve, and the irregularity is expressed in
Theorem \ref{th:q} as a linear combination of superabundances of
linear systems defined in terms of some mixed multiplier ideals
associated to this partition.  There exists a natural map from the
Galois group characters of the covering to the first orthant
appearing in the definition of the mixed multiplier ideals.  The
coefficient of each superabundance represents the number of characters
that lie in the intersection of the jumping walls associated to the
corresponding mixed multiplier ideals.  We refer the reader to
Theorem~\ref{th:q} for the precise formula and note here that 
it could be easily extended along the lines of Vaqui\'e's paper
\cite{Va}, to coverings of smooth surfaces.

The proof of Theorem \ref{th:q} occupies \S\ref{s:irregularity}.  In
\S\ref{s:applications}, the last part of the paper, some applications
are presented including E.~Hironaka's result from \cite{Hiro}
concerning the asymptotic behaviour of the irregularity of the abelian
coverings, the discussion of the general cyclic coverings, and the
computation of the irregularity of the Hirzebruch surfaces constructed
in \cite{Hi}---abelian coverings of the plane branched along
configurations of lines, \ie line arrangements.  F.~Hirzebruch mainly
deals with three arrangements and obtains three families of surfaces
with the covering group, for each family, a certain power of
$\ZZ/n\ZZ$.  For the three most interesting examples, one in each
family, namely those with $c_1^2=3c_2$, the computation of the
irregularity was performed by N.-M.~Ishida in \cite{Is}.  In
\cite{Li}, A.~Libgober computed the irregularity for two of the three
families for general $n$.  We find again one of Libgober's results,
slightly correct the second, see Proposition~\ref{p:HirzII}, and
perform the computation for the third family.

In \cite{Bu}, N.~Budur has obtained a general formula for the Hodge
numbers $h^{0,q}$, $0\leq q\leq n$, of the abelian coverings of a
smooth variety of dimension $n$.  His proof is based on the theory of
local systems of rank $1$, and the formula is expressed in terms of
the number of certain rational points inside convex polytopes (see
\cite[Theorem~1.3, Theorem~1.7]{Bu}).  A.~Libgober previously
established in \cite[\S\,3.1]{Li} a formula for the irregularity of
abelian coverings of the plane, his technique being based on mixed
Hodge structures.  The computations, mentioned above for the families
of Hirzebruch surfaces, used this formula.  His formula and ours bear
clear resemblances; it is a sum of superabundances of linear systems
expressed in terms of quasiadjunction ideals (see \cite{Li5} for the
relation between the quasiadjunction ideals and the multiplier ideals)
with coefficients given by quasiadjunction polytopes.

\section{Mixed multiplier ideals and jumping walls}
\label{s:walls}

In this section we define and characterize the jumping walls
associated to mixed multiplier ideals.  We start by briefly recalling
the notions of multiplier ideals and mixed multiplier ideals
for ideal sheaves on a smooth surface following \cite{La}. 

Let $\aaa\subseteq\Oo_X$ be a non-zero ideal sheaf on $X$ and let
$\mu:Y\to X$ be a log resolution of $\aaa$ with
$\aaa\cdot\Oo_Y=\Oo_Y(-F)$.  If $\xi$ is a positive rational number,
then the multiplier ideal associated to $\xi$ and $\aaa$ is defined as 
\[
  \Jj(\aaa^\xi) = \mu_\ast\Oo_Y(K_\mu-\floor{\xi F}).
\]
Now, for the analogous notion for several ideals, let
$\aaa_1,\ldots,\aaa_t\subset\Oo_X$ be non-zero ideals and
$\mu:Y\to X$ a common log resolution of the ideals $\aaa_i$ with
$\aaa_i\cdot\Oo_X=\Oo_Y(-F_i)$ and $\sum_iF_i+\except(\mu)$ having
simple normal crossing support.  If $\xi_1,\ldots,\xi_t$ are positive
rational numbers, then the {\it mixed multiplier ideal} associated to
the $\xi_i$ and the $\aaa_i$ is 
\[
  \Jj(\aaa_1^{\xi_1}\cdots\aaa_t^{\xi_t}) 
  = \mu_\ast\Oo_Y(K_\mu-\floor{\xi_1F_1+\cdots+\xi_tF_t}).
\]

\begin{likethm}[Definition-Lemma (see \cite{La}, Lemma 9.3.21)]
Let $\aaa\subseteq\Oo_X$ be a non-zero ideal sheaf on $X$ and let
$P\in X$ be a fixed point in the support of $\aaa$.  Then there is an
increasing sequence of positive rational numbers $\xi_j=\xi_j(\aaa,P)$
such that for every $\xi\in[\xi_j,\xi_{j+1})$,
\[
  \Jj(\aaa^{\xi_j}) = \Jj(\aaa^{\xi}) \supset \Jj(\aaa^{\xi_{j+1}}).
\]
The rational numbers $\xi_j$ are called the jumping numbers of the
ideal sheaf $\aaa$ at $P$.
\end{likethm}

The multiplier ideals and the jumping numbers are defined similarly in
the context of effective $\QQ$-divisors.  By
\cite[Proposition~9.2.28]{La}, if $C$ is a general element of the
ideal sheaf $\aaa$ and $\xi$ is a positive rational less than $1$,
then $\Jj(\xi\cdot C)=\Jj(\aaa^\xi)$.  Moreover, for any integer
divisor $C$ through a point $P$, the jumping numbers of $C$ at $P$ are
periodic and determined by the ones lying in the unit interval
$[0,1)$.  Similarly, the jumping numbers of an ideal sheaf $\aaa$ at
$P$ are periodic and determined by the ones lying in the interval
$[0,2]$.  We refer the reader to \cite[Example~9.3.24]{La} for more
ample details.

\bigskip

For the remainder of this section we consider
$\aaa_1,\ldots,\aaa_t\subset\Oo_X$ non-zero ideals such that the
subscheme defined by each $\aaa_i$ is zero dimensional and supported
at a fixed point $P\in X$.  We want to study the behaviour of the
mixed multiplier ideal $\Jj(\aaa_1^{x^1}\cdots\aaa_t^{x^t})$ as
$\xx=(x^1,\ldots,x_t)$ varies in the first orthant.  If $\mu:Y\to X$
is a log resolution defined as before, we shall denote by $E_\alpha$
the strict transforms of the exceptional divisors seen on $Y$.  There
exists effective divisors $B_\alpha$ on $Y$ such that $(B_\alpha)$ is
the dual basis to $(-E_\alpha)$ of the lattice
$\Lambda_\mu=\bigoplus_\alpha \ZZ E_\alpha$ with respect to the
intersection form on $Y$.  The basis $(B_\alpha)$ is called the {\it
branch basis} of the resolution.

Next we want to define the notion of relevant divisors.  We follow
\cite{SmTh} but see also \cite{FaJo}.

\begin{definition}
Let $\aaa\subset\mmm_P$.  A strict transform $E_\rho$ in a log
resolution of $\aaa$ is called a {\it relevant divisor} of $\aaa$ at
$P$ if either
\begin{equation} \label{eq:relevantDivisor}
  E_\rho \cdot (E_\rho^0) \geq 3,
\end{equation}
where $E_\rho^0=(\mu^\ast C)_{red}-E_\rho$ with $C$ the curve defined
by a general element of $\aaa$, or $E_\rho$ corresponds to an
arrowhead vertices of the augmented Enriques tree of $C$ at $P$.  The
index $\rho$ will be referred to as a {\it relevant position}.
\end{definition}

Note that the difference with respect
to the notion introduced in \cite{SmTh} comes from the fact that we
consider jumping numbers associated to ideal sheaves.  For example,
for the ideal of a knot, the exceptional divisor becomes a relevant
divisor.

The set of relevant positions of $\aaa$ at $P$ will be denoted by
$\RRR=\RRR_P(\aaa)$.  The following proposition stresses the
importance of the relevant divisors, or positions, in the computation
of mixed multiplier ideals.  It will further lead us to the notion of
jumping walls associated to the ideal sheaf $\aaa_1\cdots\aaa_t$ at
$P$.

\begin{pro}
  \label{pr:mainTool}
Let $\aaa_1,\ldots,\aaa_t\subset\Oo_X$ be
non-zero ideals such that the subscheme defined by each $\aaa_i$ is
zero dimensional and supported at a fixed point $P\in X$. 
Let $\mu:Y\to X$ be a log resolution of $\aaa$ and $\RRR$ the set of
relevant positions of $\aaa$ at $P$. 
It $x^i$ are positive rational numbers, then 
\[
  \Jj(\aaa_1^{x^1}\cdots\aaa_t^{x^t})
  = \mu_\ast\Oo_Y\Bigg(
  K_\mu-\sum_{\rho\in\RRR}\floor{\sum_ix^ie_i^\rho}E_\rho \Bigg),
\]
where for every $i$, 
$\aaa_i\cdot\Oo_Y=\Oo_Y(-\sum_\alpha e_i^\alpha E_\alpha)$.
\end{pro}

\proof
Consider $\yy=c\xx$ with $c\in[0,1]$.  If $c=1$ then $\yy=\xx$ and as
$c$ decreases, the coefficients $\floor{\sum_iy^ie_i^\alpha}$
decrease by discrete jumps behind.  More precisely, there is a finite
sequence of rationals $0=c_{g+1}<c_g<c_{g-1}<\cdots<c_1<c_0=1$ with
the following properties holding for any $0\leq l\leq g$:
\begin{description}\listspace
\item[\it 1)]
for any $c\in[c_{l+1},c_l)$, and for any $\alpha\not\in\RRR$,
$\floor{c_{l+1}\sum_ix^ie_i^\alpha}=\floor{c\sum_ix^ie_i^\alpha}$;
\item[\it 2)]
there exists $\BBB(l)$ disjoint from $\RRR$ such that for any
$\beta\in\BBB(l)$,  
\[
  \floor{c_{l+1}\sum_ix^ie_i^\beta}
  = \floor{c_l\sum_ix^ie_i^\beta}-1
  = c_l\sum_ix^ie_i^\beta-1;
\]
\item[\it 3)]
for any $\alpha\not\in\BBB(l)\cup\RRR$,
$\floor{c_{l+1}\sum_ix^ie_i^\alpha}=\floor{c_l\sum_ix^ie_i^\alpha}$.
\end{description} 
Set 
\[
  \Delta_l 
  = -\sum_{\alpha\not\in\RRR}\floor{c_l\sum_ix^ie_i^\alpha}E_\alpha
  -\sum_{\rho\in\RRR}\floor{\sum_ix^ie_i^\rho}E_\rho.
\]
To end the proof, it is sufficient to show that
$\mu_\ast\Oo_Y(K_\mu+\Delta_{l+1})=\mu_\ast\Oo_Y(K_\mu+\Delta_{l})$ for
any $0\leq l<g$. Set
$\Gamma=\sum_{\beta\in\BBB(l)}E_\beta$.  We have the following:

\paragraph{Claim} 
For any $\Gamma'\subset \Gamma$ and $E_\gamma\subset\Gamma'$ an
irreducible component, 
\[
  \mu_\ast\Oo_Y(K_\mu+\Delta_i+\Gamma'-E_\gamma)
  =\mu_\ast\Oo_Y(K_\mu+\Delta_i+\Gamma').
\]
We justify the claim only when $x^i$ are less than $1$.  The general
case is similar, but one needs to consider the general form of
\cite[Proposition~9.2.28]{La}.  Let $C_1,\ldots,C_t$ be the curves
defined by general elements in $\aaa_i$.  Using {\it 1)} and {\it 2)}
above we have
\[
\begin{split}
  -\Delta_l \cdot E_\gamma
  &\geq \sum_\alpha\floor{c_l\sum_i x^ie_i^\alpha}E_\alpha\cdot E_\gamma  \\
  &>
  \sum_{\beta\in\BBB(l)}c_l  \sum_ix^ie_i^\beta\, E_\beta\cdot E_\gamma
  + \sum_{\alpha\not\in\BBB(l)}
  \Big(c_l\sum_ix^ie_i^\alpha-1\Big) E_\alpha\cdot E_\gamma\\
  &\quad+ \sum_i\Big(c_lx^i-1\Big)\Ctilde_i\cdot E_\gamma\\
  &= c_l\sum_ix^i\mu^\ast C_i \cdot E_\gamma 
  - ((\mu^\ast C)_{red}-\Gamma)\cdot E_\gamma. 
\end{split}
\]
Hence 
\begin{equation}
\label{eq:le2}
  (\Delta_l+\Gamma'-E_\gamma)\cdot E_\gamma 
  < ((\mu^\ast C)_{red}-\Gamma) \cdot E_\gamma
  + (\Gamma'-E_\gamma)\cdot E_\gamma 
  \leq E_0^\gamma \cdot E_\gamma
  \leq 2
\end{equation}
since $\gamma\notin\RRR_P$.  Now, tensoring the structure sequence of
$E_\gamma$ in $Y$ with $\Oo_Y(K_\mu+\Delta_l+\Gamma')$ and pushing it
down to $X$, we get the exact sequence
\[
  0 \to
  \mu_\ast\Oo_Y(K_\mu+\Delta_l+\Gamma'-E_\gamma)
  \to
  \mu_\ast\Oo_Y(K_\mu+\Delta_{l}+\Gamma')
  \to 
  H^0(E_\gamma, K_{E_\gamma}+(\Delta_l+\Gamma'-E_\gamma)|_{E_\gamma}).
\]
The last term vanishes by (\ref{eq:le2}) justifying the claim.

From the properties {\it 2)} and {\it 3)},
$\Delta_{l+1}=\Delta_l+\Gamma$.  By repeatedly using the claim, we
obtain the result.
\qed

Next, we want to define the {\it jumping walls} associated to the
mixed multiplier ideals $\Jj(\aaa_1^{x^1}\cdots\aaa_t^{x^t})$ when
$\xx=(x^1,\ldots,x^t)$ varies in the first orthant.  The idea is that
by the previous result, such a mixed multiplier ideal varies only when
the point $\xx$ crosses certain hyperplanes defined by equations
corresponding to relevant positions.  The defining equation of such a
hyperplane is of the form
\begin{equation} \label{eq:relevantCouple}
  \sum_{i=1}^tx^ie_i^\rho = r,
\end{equation}
with $\rho\in\RRR$ and $r$ a positive integer.

\begin{definitions}
A {\it relevant value} associated to the relevant position
$\rho\in\RRR$ of the ideal $\aaa_1\cdots\aaa_t$ is a positive integer
$r$ such that there may be found a point $\yy$ in the hyperplane
$H:\sum_{i=1}^tx^ie_i^\rho=r$ and a neighbourhood $V$ of $\yy$ with
the property that the mixed multiplier ideal
$\Jj(\aaa_1^{x^1}\cdots\aaa_t^{x^t})$ corresponding to $\xx\in V$,
changes if and only if $\xx$ crosses $H$.  The pair $(\rho,r)$ is
called a {\it relevant pair} and the hyperplane $H$ a {\it jumping
wall}.
\end{definitions}

When we speak of a relevant value, we mean a positive integer which
is the relevant value associated to a certain relevant position.  Of
course, such a value might be associated to many relevant positions,
but the position we refer to will be clearly identified in the
context.

\begin{rem}
If $\aaa$ is a simple complete ideal, the relevant values
associated to the relevant position $\rho$ are the integers $\xi
e_\rho$, where $e_\rho$ is the coefficient of the strict transform
$E_\rho$ in the minimal log resolution of $\aaa$ and $\xi$ runs over
all the jumping numbers contributed by $\rho$.  We refer the reader to
\cite{Ja,Na2} for a formula producing all these jumping numbers.
\end{rem}

\begin{exa}
Let $\aaa_1=(u^3,v^2)$ and $\aaa_2=(u^6,v^2)$ be ideals in $\CC[u,v]$.
Let $E_1,E_2$ and $E_3$ be the exceptional divisors necessary for the
minimal log resolution of $\aaa_1$ and let $E_4$ be the supplementary
exceptional divisor necessary for finishing the minimal log resolution
of $\aaa_2$.  Clearly, if $C_i$ are general elements in each $\aaa_i$,
then $\mu^\ast(C_1+C_2)=\Ctilde_1+\Ctilde_2+4E_1+7E_2+12E_3+9E_4$.
The divisors $E_3$ and $E_4$ are the only relevant divisors.  Then $5$
and $7$ are the first relevant values associated to the relevant
divisor $E_3$ with the jumping walls $H_{(3,r)}:6x^1+6x^2=r$, $r=5,7$.
Moreover, $4$ and $5$ are the first relevant values associated to
$E_4$ with the jumping walls $H_{(4,r)}:3x^1+6x^2=r$, $r=4,5$.

The point $\yy$ from the definition of the jumping wall $H_{(4,4)}$
can by any point on $H_{(4,4)}\cap \RR_+^2$ with $y^1<1/3$.  The other
points in the intersection do not satisfy the property in the
definition of the relevant value.  If $y^1>1/3$ then on a sufficiently
small neighbourhood of $\yy$, the mixed multiplier ideal
$\Jj(x^1C_1+x^2C_2)$ equals the maximal ideal $(u,v)$.  If
$\yy=(1/3,2/3)$ then the multiplier ideal also changes when it crosses
the wall $H_{(3,5)}:6x^1+6x^2=5$.  In the figure above, if $\xx$ lies
in the open shaded polygon, then the mixed multiplier ideal equals the
maximal ideal.
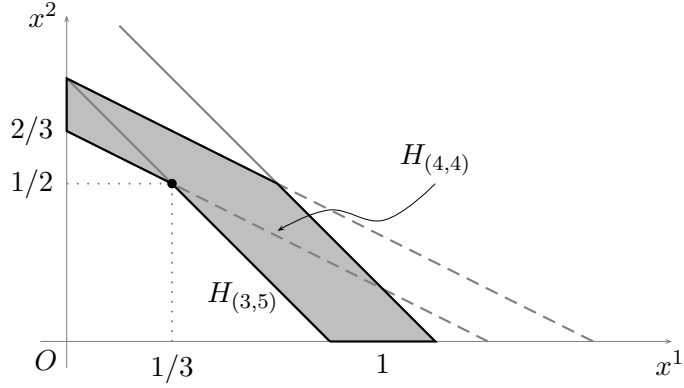
\begin{center}
  \begin{pspicture}(-1,-.5)(7.5,4.7) 
    \psset{unit=.7cm}
    \pscustom[fillstyle=solid,fillcolor=lightgray,linecolor=lightgray]{%
      \psline(5,0)(7,0)(4,3)(0,5)(0,4)(2,3)}

    \psline[linewidth=0.25pt,linecolor=gray,arrows=->](-.5,0)(11.5,0)
    \psline[linewidth=0.25pt,linecolor=gray,arrows=->](0,-.5)(0,6.2)
    \psline[linestyle=dotted,linecolor=gray](2,3)(2,0)
    \psline[linestyle=dotted,linecolor=gray](2,3)(0,3)
    \psline[linecolor=gray](5,0)(0,5)
    \psline[linecolor=gray](7,0)(1,6)
    \psline[linecolor=gray,linestyle=dashed](8,0)(0,4)
    \psline[linecolor=gray,linestyle=dashed](10,0)(0,5)
    \psline(5,0)(7,0)(4,3)(0,5)(0,4)(2,3)(5,0)
    \psdot(2,3)

    \uput{\labelsep}[dl](4.2,1.3){$H_{(3,5)}$}
    \uput{\labelsep}[u](7,2.9){$H_{(4,4)}$}
    \pscurve[linewidth=0.25pt]{<-}(4,2.1)(5,2.5)(6,2.3)(7,3)
    \uput{\labelsep}[dl](0,0){$O$}
    \uput{\labelsep}[dr](11,0){$x^1$}
    \uput{\labelsep}[ul](0,5.8){$x^2$}
    \uput{\labelsep}[d](2,0){$1/3$}
    \uput{\labelsep}[d](6,0){$1$}
    \uput{\labelsep}[l](0,4){$2/3$}
    \uput{\labelsep}[l](0,3){$1/2$}
  \end{pspicture}
\end{center}
\end{exa}

For practical reasons, what we have to do next is to determine a
relatively small set of candidates for the relevant values associated
to $\rho$.

\begin{definition}
  \label{d:relIdeals}
Let $\rho$ be a relevant position for the ideal $\aaa$ and $\mu$ a log
resolution.  The relevant ideal associated to $\aaa$ and $\rho$
is the simple complete ideal $\mu_\ast\Oo_Y(-B_\rho)$, where $B_\rho$
is the $rho$ element in the branch basis of the resolution.
\end{definition}

\begin{pro}
  \label{pr:mainR}
Let $\aaa_1,\ldots,\aaa_t\subset\Oo_X$ be non-zero ideals such that
the subscheme defined by each $\aaa_i$ is zero dimensional and
supported at a fixed point $P\in X$.  Let $\mu:Y\to X$ be a log
resolution of $\aaa=\aaa_1\cdots\aaa_t$ with $(B_\alpha)$ the branch
basis of the resolution.  Then the set of relevant values associated
to the relevant position $\rho$ is contained in the set of
relevant values associated to $\rho$ of the relevant ideal 
$\mu_\ast\Oo_Y(-B_\rho)$.
\end{pro}

\proof
It is sufficient to consider the case $t\geq2$.  Let $\rho_0$ be a
relevant position and $r$ a relevant value with
$H:\sum_{i=1}^tx^ie_i^{\rho_0}=r$ the corresponding hyperplane.  The
point $\yy$ may be chosen such that $H$ is the only jumping hyperplane
containing it.  It is here that we need $t\geq 2$.  Using
Proposition~\ref{pr:mainTool}, since 
\[
  \mu_\ast\Oo_Y(K_\mu-\floor{\sum_iy^iF_i}) 
  \subset \mu_\ast\Oo_Y(K_\mu-\floor{\sum_iy^iF_i}+E_{\rho_0})
  = \mu_\ast\Oo_Y(K_\mu-\floor{\sum_i(1-\epsilon)y^iF_i})
\]
with $0<\epsilon\ll 1$, we get that
\[
  \mu_\ast\Oo_Y(K_\mu-\sum_{\rho\in\RRR}r^\rho E_\rho)
  \subset \mu_\ast\Oo_Y(K_\mu-\sum_{\rho\neq\rho_0}r^\rho E_\rho
    -(r-1)E_{\rho_0}).
\]
Setting $K_\mu=\sum_\alpha k^\alpha E_\alpha$ and
$\bbb=\mu_\ast\Oo_Y(-\sum_{\rho\neq\rho_0}(r^\rho-k^\rho)E_\rho)$, it
follows that
\[
  \bbb\cap\mu_\ast\Oo_Y((k^{\rho_0}-r)E_{\rho_0})
  \subset \bbb\cap\mu_\ast\Oo_Y((k^{\rho_0}-r+1)E_{\rho_0}),
\]
\ie that 
$\mu_\ast\Oo_Y((k^{\rho_0}-r)E_{\rho_0})\subset
\mu_\ast\Oo_Y((k^{\rho_0}-r+1)E_{\rho_0})$.  
Now, set $\qqq=\mu_\ast\Oo_Y(-B_{\rho_0})$.  The Enriques tree
associated to $\mu':Y'\to X$, the minimal log resolution of $\qqq$, is
the path from the root to the vertex $P_{\rho_0}$ of the Enriques tree
associated to $\aaa$.  Let $\VVV'$ be the set of vertices of this path
and $\RRR'\subset\RRR\cap\VVV'$ the set of relevant positions.  If
$\qqq\cdot\Oo_{Y'}=\Oo_{Y'}(-\sum_{\alpha\in\VVV'}e^\alpha E_\alpha)$,
let $\RRR''\subset\RRR'$ be the subset of relevant positions such that
for any $\rho\in\RRR''$, $re^{\rho}/e^{\rho_0}$ is an integer.  Then,
using again Proposition~\ref{pr:mainTool} and the previous strict
inclusion,
\begin{multline*}
  \Jj(\qqq^{r/e^{\rho_0}})
  = \mu_\ast\Oo_{Y'}\Bigg(
    K_{\mu'}-\sum_{\rho\in\RRR'}\frac{re^\rho}{e^{\rho_0}}E_\rho\Bigg)
  = \mu_\ast\Oo_{Y'}\Bigg(
    \sum_{\rho\in\RRR'}\bigg(k^\rho-\frac{re^\rho}{e^{\rho_0}}\bigg)
    E_\rho\Bigg) \\
  \subset \mu_\ast\Oo_{Y'}\Bigg(
    \sum_{\rho\in\RRR'\smallsetminus\RRR''}
    \bigg(k^\rho-\frac{re^\rho}{e^{\rho_0}}\bigg)E_\rho+
    \sum_{\rho\in\RRR''}
    \bigg(k^\rho-\frac{re^\rho}{e^{\rho_0}}+1\bigg)
    E_\rho\Bigg)
  = \Jj(\qqq^{(1-\epsilon)r/e^{\rho_0}}),
\end{multline*}
where $0<\epsilon\ll 1$.  Hence $r/e^{\rho_0}$ is a jumping number of
$\mu_\ast\Oo_Y(-B_{\rho_0})$ associated to the relevant position
$\rho_0$. 
\qed

\begin{exa}
Let $\aaa_1$ and $\aaa_2$ be simple complete ideals supported at $P$.
We want to show that for a relevant position $\rho$, the set of
relevant values of the ideal $\mu_\ast\Oo_Y(-B_\rho)$ are indeed
needed, \ie that the union of the sets of of relevant values of each
ideal $\aaa_i$ associated to $\rho$ is not sufficient.  Let $C_1$ and
$C_2$ be two unibranch curves, general elements in $\aaa_1$ and
$\aaa_2$, and let the associated augmented Enriques tree of the
minimal log resolution of $C_1+C_2$ be as in the figure below.
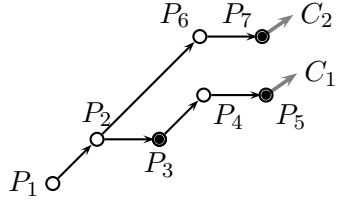
\begin{center}
\begin{pspicture}(0,.2)(4,2.8)
  \psset{arrows=->,radius=.1,unit=5ex}
  \rput(1,1){
    \psline[linewidth=.5mm,linecolor=gray]{->}(2.65,1.65)(3.15,2)

    \Cnode(-.71,-.71){1}
    \Cnode(0,0){2}
    \Cnode(1,0){3}
    \dotnode(1,0){3}
    \Cnode(1.65,1.65){4}
    \Cnode(2.65,1.65){5}
    \dotnode(2.65,1.65){5}
    
    \ncline{1}{2}
    \ncline{2}{3}
    \ncline{2}{4}
    \ncline{4}{5}
    \uput{\labelsep}[l](-.71,-.71){$P_1$}
    \uput{\labelsep}[u](0,0){$P_2$}
    \uput{\labelsep}[d](1,0){$P_3$}
    \uput{\labelsep}[ul](1.65,1.65){$P_6$}
    \uput{\labelsep}[ul](2.65,1.65){$P_7$}
    \uput{\labelsep}[r](3.05,2){$C_2$}
  }
  \rput(2.71,1.71){
    \psline[linewidth=.5mm,linecolor=gray]{->}(1,0)(1.50,.35)

    \Cnode(-.71,-.71){1}
    \Cnode(0,0){2}
    \Cnode(1,0){3}
    \dotnode(1,0){3}
    \ncline{1}{2}
    \ncline{2}{3}
    \uput{\labelsep}[dr](0,0){$P_4$}
    \uput{\labelsep}[dr](1,0){$P_5$}
    \uput{\labelsep}[r](1.40,.35){$C_1$}}
\end{pspicture} 
\end{center}
We have 
$\mu^\ast(C_1+C_2)=\Ctilde_1+\Ctilde_2
+6E_1+10E_2+18E_3+19E_4+38E_5+11E_6+22E_7$.  
The relevant positions are indicated by the black vertices: $3$, $5$
and $7$.  The jumping numbers of $\aaa_1$ contributed by $E_3$ are
$(5+6k)/12$, with $k\in\NN$.  But the three first jumping numbers of
$\aaa_1\aaa_2$ are $5/18$, $7/18$ and $8/18$.  Hence the relevant
values $7$ and $8$ are not among the relevant values of $\aaa_1$
associated to the relevant position $3$.  Of course, the well known
jumping numbers of the ideal $\mu_\ast\Oo_Y(-B_3)$ are $(2a+3b)/6$,
with $a$ and $b$ positive integers.
\end{exa}

\section{The irregularity of the abelian covering 
of the projective plane}
\label{s:irregularity}

In this section we state and prove in Theorem~\ref{th:q} a formula for
the irregularity of the standard abelian coverings of the plane.  To
be able to express it, we start by summarizing the definition and some
properties of these coverings in a form convenient for further use.
Then, using the jumping walls in the context of plane curves, we
introduce the notion of distinguished faces, leading notion in formula
(\ref{eq:q}).

\subsection{Abelian coverings}
\label{ss:abelianCoverings}
Let $\pi:Y\to X=Y/G$ be a Galois covering with abelian Galois group
$G$.  It is well known that $\pi_\ast\Oo_Y$ is a coherent sheaf of
$\Oo_X$-algebras and that $Y\simeq\Specbf_{\Oo_X}(\pi_\ast\Oo_Y)$.  In
addition, if $Y$ is normal and $X$ is smooth, $\pi$ is flat and
consequently $\pi_\ast\Oo_Y$ is locally free of rank $n$.  The action
of $G$ on $\pi_\ast\Oo_Y$ decomposes it into the direct sum of eigen
line bundles associated to the characters
$\chi\in\Ghat=\Hom(G,\SS^1)$,
\begin{equation*}
  \pi_\ast\Oo_Y
 = \Oo_X\oplus\bigoplus_{\chi\in\Ghat,\chi\neq1}\Ll_\chi\inv.
\end{equation*}
The action of $G$ on $\Ll_\chi\inv$ is the multiplication by
$\chi$.

Let $\chi_1,\ldots,\chi_s\in\Ghat$ such that the group of characters
is the direct sum of the cyclic subgroups generated by
$\chi_1,\ldots,\chi_s$.  Let $n_1,\ldots,n_s$ be their orders.  In
\cite[Proposition 2.1]{Pa} it is shown that the ring structure of
$\pi_\ast\Oo_Y$, and hence $Y$, are determined by the following linear
equivalences or isomorphisms.  For every $1\leq j\leq s$,
\begin{equation}
  \label{eq:RBDRel}
  n_jL_{\chi_j} 
  \sim \sum_{f\in\FFF}\frac{n_jf(\chi_j)^\bullet}{m_f}\,B_f,
\end{equation}
where: 1) the set $\FFF$ consists of all group epimorphisms from
$\Ghat$ to different $\ZZ/m\ZZ$; 2) the curve $B_f\subset X$ with
$f\in\FFF$ is the sub-divisor of the branch locus defined
set-theoretically as $\pi(R_f)$, with $R_f$ the union of all the
components $D$ of the ramification locus associated to the group
epimorphism $f$; 3) the integer $a^\bullet$ denotes the smallest
non-negative integer in the equivalence class $a\in\ZZ/m$ (each time
the integer $m$ being understood from the context).

\begin{rem*}
If $D$ is a component of the
ramification locus, since $Y$ is normal and $X$ smooth, $D$ is
$1$-codimensional.  The {\it inertia subgroup} $H\subset G$ and a
character $\psi\in\Hhat$---the induced representation of $H$ on the
cotangent space to $Y$ at $D$---that generates $\Hhat$ are associated
to $D$.  Dualizing the inclusion $H\subset G$, such a pair $(H,\psi)$
is equivalent to a group epimorphism $f\colon\Ghat\to\ZZ/m_f$, where
$m_f=|H|$. 
\end{rem*}

Following \cite{Pa}, the line bundles $\Ll_{\chi_j}$, $1\leq j\leq s$,
and the divisors $B_f$, $f\in\FFF$, are called a set of {\it reduced
building data} for the covering.  In case $X$ is compact, the covering
is uniquely determined by the isomorphisms (\ref{eq:RBDRel}), up to
isomorphisms of abelian coverings.  It is to be noticed that if
$\chi=\chi_1^{a_1}\cdots\chi_s^{a_s}\in\Ghat$, then
\begin{equation} \label{eq:L_chi}
  L_\chi 
  \sim \sum_{j=1}^s a_jL_{\chi_j}
  -\sum_{f\in\FFF}
  \floor{\sum_{j=1}^s\frac{a_jf(\chi_j)^\bullet}{m_f}}B_f.
\end{equation}
So, to a normal covering $\pi:Y\to X$ of a smooth variety $X$, a set
of reduced building data is associated satisfying the relations
(\ref{eq:RBDRel}).  Conversely, starting with a set of reduced
building data and relations (\ref{eq:RBDRel}) an abelian covering is
constructed which will be called a {\it standard abelian covering}.

\begin{notation}
  \label{n:theSurface}
In the sequel, we fix the notation $S(\nn,M,\cC,H_\infty)$ for the
abelian covering of the plane constructed as follows.  Let
$C\subset\PP^2$ be a reduced curve and $H_\infty\subset\PP^2$ a line
called the line at infinity.  The set of reduced building data
consists of
\begin{itemize}\listspace
\item
the line bundles 
$\Ll_{\chi_j}=\Oo_{\PP^2}(\ceil{\sum_{i=1}^t\mu_j^id_i/n_j})$, 
$1\leq j\leq s$,
\item
the line $H_\infty$ and the curves $C_i\subset C$ of degree $d_i$,
$1\leq i\leq t$, such that $C=\sum_iC_i$,
\item
the linear equivalences 
$n_jL_{\chi_j} \sim
\sum_{i=1}^t\mu_j^iC_i+(\ceil{\mu_jd_j/n_j}n_j-\mu_jd_j)H_\infty$,
$1\leq j\leq s$.
\end{itemize}
The covering $S\to\PP^2$ thus obtained has Galois group
$\oplus_{j=1}^s\ZZ/n_j\ZZ$ and depends on the line at infinity and the
$r\times s$ matrix $M=[\mu_j^i]$ with non negative integer entries.
$\cC$ is the list of curves $(C_1,\ldots,C_t)$ and $\nn$ the
$s$-vector $(n_1,\ldots,n_s)$ defining the covering group.  The list
of curves $\cC$ such that $C=\sum_iC_i$ will be referred to as a {\it
partition} of $C$.
\end{notation}

\begin{rem}
  \label{r:hironaka}
In \cite{Hiro} the construction of the abelian coverings $\Sigma_n$
that are studied is different from the one presented above.  
The construction is based on the composition of the Hurewich
epimorphism with the morphism given by the change of constants in the
homology groups
\[
  \pi_1(U) \lra H_1(U,\ZZ) \lra H_1(U,\ZZ/n\ZZ).
\]
Here $U=\PP^2\smallsetminus(C\cup H_\infty)=\CC^2\smallsetminus C$,
with $C=\sum_{j=1}^sC_j$ the decomposition of $C$ into {\it
irreducible components}.  It is known that there exists an exact
sequence (see \cite[Proposition~1.3]{Di})
\[
  \ZZ \stackrel{\iota}{\lra} \bigoplus_{j=0}^s \ZZ 
  \lra H_1(U,\ZZ) \lra 0
\] 
with $\iota(1)=g_0+\sum_{j=1}^sd_jg_j$, $g_0=(1,0,\ldots)$ and so on.
It follows that the epimorphism corresponds to a Galois unbranched
covering $V\to U$ with group $H_1(U,\ZZ/n\ZZ)\simeq(\ZZ/n\ZZ)^s$.  By
the existence theorem of Grauert and Remmert \cite{GrRe}, this
covering extends to a unique normal abelian covering
$\pi:\Sigma_n\to\PP^2$, \ie such that $\pi\inv(U)=V$.  It turns out
that $\Sigma_n$ is the normalization of the standard covering
$S((n,\ldots,n),I_s,\cC,H_\infty)$ introduced in
Notation~\ref{n:theSurface}, \ie for which the linear equivalences
(\ref{eq:RBDRel}) are given by 
$nL_{\chi_j}\sim C_j+(\ceil{d_j/n}n-d_j)H_\infty$.
\end{rem}

\subsection{Jumping walls and distinguished faces 
of a curve endowed with a partition}

Let $C$ be a reduced plane curve endowed with a partition
$\cC=(C_1,\ldots,C_t)$.  Note that $C=\sum_iC_i$.  For each singular
point of $C$, we shall need to consider the mixed multiplier ideal
$\Jj(\xx\cdot \cC)=\Jj(x^1C_1+\ldots x^tC_t)$ with $\xx$ varying in
the hypercube $[0,1)^t$.  Interpreting the results of \S~\ref{s:walls}
in this context, it follows that the jumping walls associated to these
multiplier ideals cut up the hypercube into convex rational polytopes
on which the map $x\mapsto\Jj(\xx\cdot\cC)$ is constant.  Note that
the fibres of this map are neither open nor closed.

\begin{definition}
A {\it face}\footnote{In a previous
version of the paper these faces were called walls and A.~Libgober
kindly pointed to me that using the word {\it wall} was misleading in
codimension $\geq 2$.} associated to $C$ endowed with the partition
$\cC$ is a finite intersection of jumping walls and coordinate
hyperplanes.
\end{definition}

\begin{rem}
Even if we are interested in the jumping walls intersecting the
hypercube, the context of rational divisors is not sufficient, since
the jumping walls are determined by relevant values associated to
ideal sheaves---more precisely, a jumping wall might intersect a
coordinate axis in a jumping number bigger than $1$ associated only to
an ideal sheaf.
\end{rem}

If $W$ is a face associated to a curve $C$ endowed with the partition
$\cC$, the set $\Uu(W)$ is the set of the connected components of the
difference between $W$ and the union of all the jumping walls and
coordinate hyperplanes that do not contain $W$.  The mixed multiplier
ideal is constant on each $U\in\Uu(W)$ and will be denoted by
$\Jj(U\cdot\cC)$.  Furthermore, if $\dd$ is the vector
$(d_1,\ldots,d_t)$, where $\deg C_i=d_i$, we define the height
function $h_\cC:\RR^t\to\RR$ by $h_\cC(\xx)=\dd\cdot\xx$.

\begin{definition}
A face $W$ of the projective curve $C$ endowed with the partition 
$\cC$ is called a {\it distinguished face} if the height function
$h_\cC$ is constant on $W$.  The set of distinguished faces 
will be denoted by $\FFF(\cC)$.
\end{definition}

\begin{lem}
If $C\subset\PP^2$ is endowed with the partition $\cC$, $C_i$ is a 
component in $\cC$ and $W$ a distinguished face, then either $x^i=0$ 
along $W$, or $C_i$ passes through $P$, a singular point of $C$, 
to which one of the jumping walls that cut out $W$ is associated.
\end{lem}

\proof
It is easy to see that a distinguished face $W$, seen as a subset in
the first orthant of $\RR^t$, is bounded for the euclidean metric.
Suppose that $W\not\subset\{x^i=0\}$.  Then the component $C_i$ must
satisfy the conclusion since otherwise the corresponding coordinate
$x^i$ would be unbounded. 
\qed

\begin{exa}
  \label{ex:ceva}
Let $C_1,\ldots,C_6$ be the lines of Ceva's arrangement $C=\sum C_i$;
$C_i$ and $C_j$ intersect in a node of the arrangement if and only if
$i+j=7$.  Ceva's arrangement has four triple points: 
$C_4\cap C_5\cap C_6$, $C_1\cap C_2\cap C_4$, $C_2\cap C_3\cap C_6$
and $C_1\cap C_3\cap C_5$.  For $C$ endowed with the partition
$\cC=(C_1,\ldots,C_6)$ there are five distinguished faces: one for
each triple point and one for all four.  Clearly for each triple point
$P$ there is a distinguished face $W_P$; for example if 
$P=C_4\cap C_5\cap C_6$ then $W_P$ is defined by $x^4+x^5+x^6=2$,
$x^1=x^2=x^3=0$ and $h_\cC(W_P)=2$.  Now, if $W$ is a distinguished
face different from the $W_P$, then let $\phi_\alpha(\xx)=2$ be the
equations defining the jumping walls that cut out $W$---$2$ is the
only relevant value.  Note that each equation is of the form
$x^i+x^j+x^k=2$.  Let $I\subset\{1,2,\ldots,6\}$ be the set of
subscripts appearing in the equations $\phi_\alpha$.  Since $W$ is
distinguished, $x^j=0$ along $W$ for every $j\not\in I$.  Furthermore
the equation
\[
  \sum_{i\in I} x^i = h_\cC(W)
\]
is a linear combination of the $\phi_\alpha$.  Hence there exist
$\zeta_\alpha$ such that 
\[
  \sum_\alpha \zeta_\alpha (\phi_\alpha(\xx)-2) 
  = \sum_{i\in I} x^i - h_\cC(W)
\]
for any $\xx\in\RR^6$.  Hence $2\sum_\alpha\zeta_\alpha=h_\cC(W)$, and 
taking $x^i=1$ for every $i\in I$, 
$3\sum_\alpha\zeta_\alpha = |I|$.  It follows that $2|I|=3h_\cC(W)$, 
\ie that $|I|=6$ and $h_\cC(W)=4$.  To see that $W$ is unique with these
properties it is sufficient to notice that $W$ is defined by the four 
equations corresponding to the four triple points.  It is clear that it 
should be defined by at least three out of four equations.  Summing these
three equations and using $\sum_1^6x^i=4$ we get the fourth.
\end{exa}

\begin{exa}
  \label{exa:2tacNodes}
Let $\Gamma_1$ and $\Gamma_2$ be two conics that have common tangents
at the two points of intersection $P$ and $Q$.  Let $H_\infty$ be the
line through $P$ and $Q$.  We want to determine the set of
distinguished faces $\FFF_d$ for the curve $C=C_1+C_2$, with the
partition $\cC=\{C_1,C_2\}$, where $C_1=\Gamma_1+\Gamma_2$ and
$C_2=H_\infty$.  The curve $C_1$ has two tacnodes at $P$ and $Q$; a
jumping number $3/4$ and hence a unique relevant value $3$.  The curve
$C=C_1+C_2$ has two singular points and the exceptional configuration
of the minimal log-resolution is $(2+1)E_1+(4+1)E_2$.  There are two
jumping values, $3/5$ and $4/5$ and two relevant values $3$ and $4$
associated to the second exceptional divisor in the log-resolution for
each singular point.  It follows that there are two jumping walls
$W_3$ and $W_4$ defined by $4x^1+x^2=3$ and $4x^1+x^2=4$ respectively.
There are three faces and all three are distinguished since
$h_\cC=4x^1+x^2$: $W_3$, $W_4$ and the intersection of $W_3$ with the
coordinate line $\{x^2=0\}$, \ie the point $W_0$ of coordinates
$(3/4,0)$.  Finally,
\[
  \Uu(W_0)=\{W_0\},\quad
  \Uu(W_3)=\{W_3\smallsetminus W_0\}
  \textq{and}
  \Uu(W_4)=\{W_4\}.
\]
\end{exa}

\subsection{The irregularity}

In this section we state and prove the formula for the irregularity of
the abelian covering $S'=S(\nn,M,\cC,H_\infty)$---the standard
$\oplus_{j=1}^s\ZZ/n_j\ZZ$-covering $S'\to\PP^2$ defined by the linear
equivalences
\[
  n_jL_{\chi_j} \sim \sum_{i=1}^t\mu_j^iC_i+
  \bigg(
  \ceil{\frac{1}{n_j}\sum_{i=1}^t\mu_j^id_i}n_j-\sum_{i=1}^t\mu_j^id_i
  \bigg)H_\infty,
\]
where $\Ll_{\chi_j}=\Oo_{\PP^2}(\ceil{\sum_{i=1}^t\mu_j^id_i/n_j})$,
$d_i=\deg C_i$, $\cC=(C_1,\ldots,C_t)$, $\nn=(n_1,\ldots,n_s)$ and $M$
denotes the $t\times s$ matrix $[\mu_j^i]$ (see
Notation~\ref{n:theSurface}).

For any rational convex polytope $U\subset \RR^t$ set  
\[
  |U|_\nn^M = \card\phi\inv(U\cap[0,1)^t),
\]
where the map
$\phi : [0,1)^s\cap\bigoplus_{j=1}^s1/n_j\ZZ \lra [0,1)^t$, depending  
on $\nn=(n_1,\ldots,n_s)$ and the matrix $M$, is defined by
\[
  \phi\bigg(\frac{a^1}{n_1},\ldots,\frac{a^s}{n_s}\bigg) 
  = \bigg(\!\fpart{\sum_j\mu_j^1\frac{a^j}{n_j}},\ldots,
    \fpart{\sum_j\mu_j^t\frac{a^j}{n_j}}\!\bigg).
\]
In case $M$ is the identity matrix we shall omit the superscript $M$
it the notation $|U|_\nn^M$.  Similarly, if $n_j=n$ for every $j$, we
shall use $|U|_n^M$ for $|U|_{(n,\ldots,n)}^M$.  Note that $|U|_n$ is
the number of rational points in $W$, points whose coordinates belong
to $1/n\ZZ$ and are non-negative and smaller than $1$.

\begin{thm} 
  \label{th:q}
Let $S$ be the normalization of the $\oplus_{j=1}^s\ZZ/n_j\ZZ$-abelian
covering $S'\to\PP^2$ with $S'=S(\nn,M,\cC,H_\infty)$.  Suppose that
$S$ is a connected surface. Then
\begin{equation} \label{eq:q}
  q(S) = 
  \sum_{W\in\FFF(\cC)}\sum_{U\in\Uu(W)}|U|_\nn^M
  \cdot h^1\big(\PP^2,\Oo_{\PP^2}(-3+h_{\cC}(W))
  \otimes\Jj(U\cdot\cC)\big),
\end{equation}
if $\sum_iC_i$ is transverse to $H_\infty$, and 
\begin{equation} \label{eq:qbar}
  q(S) = 
  \sum_{W\in\FFF(\overline{\cC})}\sum_{U\in\Uu(W)}|U|_\nn^\Mbar
  \cdot h^1\big(\PP^2,\Oo_{\PP^2}(-3+h_{\overline{\cC}}(W))
  \otimes\Jj(U\cdot\overline{\cC})\big),
\end{equation}
if $\sum_iC_i$ is not transverse to $H_\infty$ and the covering is
branched along $H_\infty$, where $\overline{\cC}=(H_\infty,\cC)$ and
\[
  \Mbar = 
  \begin{bmatrix}
    \mu_1^0 & \cdots & \mu_s^0 \\
    & M 
  \end{bmatrix}
\]
with $\mu_j^0=\ceil{\sum_{i=1}^t\mu_j^id_i/n_j}n_j-\sum_{i=1}^t\mu_j^id_i$.
\end{thm}

\proof
In order to compute the irregularity of $S$ we need to see $S$ as a
standard abelian covering of group $\oplus_{j=1}^s\ZZ/n_j\ZZ$.  We use
the normalization algorithm from \cite{Pa}.  Let $\mu:X\to\PP^2$ be a
log resolution of the branch divisor.  According to the position of
the line at infinity, the points that are blown up lie either on
$\sum_iC_i$ or on $\sum_iC_i+H_\infty$.  The abelian covering
$S'\to\PP^2$ pulls back to a standard abelian covering $S''\to X$
defined by line bundles $\Ll''_{\chi_j}$.  Then, the normalization
procedure yields the normal surface $S$ with only Hirzebruch-Jung
singularities.
\[
\begin{CD}
  S @>>> S'' @>>> S' \\
  @VV{\pi}V   @VVV     @VVV \\
  X @= X @>\mu>> \PP^2 \\
\end{CD}
\]
It is 
a standard abelian covering with line bundles $\Ll_{\chi_j}$
among the elements of the reduced 
building data.  Using the Leray spectral sequence
and the Serre duality,
\[
  q(S) 
  = h^1(S,\Oo_S)
  = h^1(X,\pi_\ast\Oo_S)
  = \sum_{\chi\in\Ghat}h^1(X,\omega_X\otimes\Ll_\chi).
\]
For the computation of the terms in the right hand member, we
distinguish two cases.

\paragraph{First case} $H_\infty$ is transverse to $\sum_iC_i$.
Let  
\begin{equation}
  \label{eq:pull-back}
  \mu^\ast C_i = \Ctilde_i+\sum_P\ee_i^P\cdot\eE_P
\end{equation}
the sum being taken over all the singular points of $\sum_iC_i$
excepting the nodes.
Here and in the sequel $\ee_i^P\cdot\eE_P$ denotes the sum
\[
  \sum_\alpha e_i^{P,\alpha} E_{P,\alpha}
\]
where $E_{P,\alpha}$ are the irreducible components of the exceptional
configuration of the log resolution $\mu$ over $P$.  The line bundles
$L''_{\chi_j}\sim \ceil{\sum_i\mu_j^id_i/n_j}\Htilde$ and the linear
equivalences
\[
  n_jL''_{\chi_j}
  \sim  \sum_i\mu_j^i\Ctilde_i + 
  \sum_{i,P}\mu_j^i\ee_i^P\cdot\eE_{P}
  + \Bigg(\!\ceil{\sum_i\mu_j^id_i/n_j}n_j-\sum_i\mu_j^id_i\Bigg)
  \Htilde_\infty,
\]
holding for each $1\leq j\leq s$,
define $S''$.  If $\chi=\chi_1^{a^1}\cdots\chi_s^{a^s}$,
$0\leq a^j<n_j$, then, after normalization, by 
\cite[Proposition 3.2]{Na} and by Proposition~\ref{p:normalization4},
$L_\chi$ is linearly equivalent to
\begin{multline*}
  \sum_j a^jL''_{\chi_j}
  -\sum_i\floor{\sum_j\frac{a^j\mu_j^i}{n_j}}\Ctilde_i
  -\sum_P\floor{\sum_{i,j}\frac{a^j\mu_j^i}{n_j}\ee_i^P}\cdot\eE_P \\
  -\floor{\sum_{j=1}^s\frac{a^j}{n_j}
  \Big(\!\ceil{\sum_i\mu_j^id_i/n_j}n_j-\sum_i\mu_j^id_i\Big)}
    \Htilde_\infty,
\end{multline*}
and using (\ref{eq:pull-back}) and $\Htilde_\infty\sim\Htilde$ to
\[
  \Bigg(\!
  \ceil{\sum_{i,j}\frac{a^j\mu_j^id_i}{n_j}}
  -\sum_i\floor{\sum_j\frac{a^j\mu_j^i}{n_j}}d_i
  \Bigg)\Htilde 
  -\Bigg(\!
  \sum_P\floor{\sum_{i,j}\frac{a^j\mu_j^i}{n_j}\ee_i^P}
  +\sum_{P,i}\floor{\sum_j\frac{a^j\mu_j^i}{n_j}}\ee_i^P
  \Bigg)\cdot\eE_P.
\]
Setting 
\begin{equation}
  \label{eq:theX}
  x^i = \fpart{\sum_ja^j\mu_j^i/n_j}
\end{equation}
for every $1\leq i\leq t$, it follows that 
\[
  L_\chi \sim \ceil{\sum_id_ix^i}\Htilde 
  -\sum_P\floor{\sum_ix^i\ee_i^P}\cdot\eE_P.
\]
Then
\begin{equation} \label{eq:local}
  h^1(X,K_X+L_\chi) = 
  h^1(\PP^2,\Oo_{\PP^2}(-3+\ceil{\sum_id_ix^i})\otimes\Jj(\sum_ix^iC_i))
\end{equation}
and the dimension $h^1(X,K_X+L_\chi)$ might be non-zero whenever the
numbers $x^i$ satisfy three conditions.  First $\sum_ix^i d_i$ must be
an integer.  If not, then the right-hand side of (\ref{eq:local})
vanishes by the Kawamata-Viehweg-Nadel vanishing theorem.  Second, for
every curve $C_i$ there exists a singular point $P$ of $B$ lying on
$C_i$ and a relevant position $\alpha$ of $P$ such that the number
$\sum_ix^i e_i^{P,\alpha}$ is a relevant value of $D$ at $(P,\alpha)$.
Indeed, if this condition does not hold for $C_1$ for example, it is
sufficient to notice that $\Jj(D')=\Jj(\sum_ix^iC_i)$, where
$D'=(x^1-\epsilon)C_1+\sum_{i=2}^tx^i C_i$ for $\epsilon>0$
sufficiently small, and to apply the Kawamata-Viehweg-Nadel vanishing
theorem to see that the right-hand side of (\ref{eq:local}) vanishes.
Third, suppose that $\sum_ix^id_i$ is an integer and that for every
component $C_i$ there exists a singular point $P$ of $B$ and a
position $\alpha$ such that $r_{P,\alpha}=\sum_ix^ie_i^{P,\alpha}$ is
a relevant value of $D$ at $(P,\alpha)$.  If $W$ is the space of
solutions of these equations seen as equations in the unknowns $x^i$,
then $W$ is a face for the partition $\cC=(C_1,\ldots,C_t)$ and the
linear operator $h_\cC:\xx\mapsto\sum_ix^id_i$ must be constant on
$W$.  Indeed, if $W$ is positive dimensional and not contained into a
fibre of $h_\cC$, it is sufficient to take $\yy\in W$ such that
$\sum_iy^id_i<\sum_ix^id_i$.  Then, if
$\delta=\ceil{\sum_iy^id_i}=\sum_ix^id_i$,
\[
  h^1(\PP^2,\Oo_X(-3+\sum_id_ix^i)\otimes\Jj(\sum_ix^iC_i)) =
  h^1(\PP^2,\Oo_X(-3+\delta)\otimes\Jj(\sum_iy^iC_i)) =
  0.
\]
So the face $W$ is distinguished and
$\Jj(\sum_ix^iC_i)=\Jj(U\cdot\cC)$, for $U$ the corresponding
connected components defined on $W$ by the other jumping walls and
coordinate hyperplanes.

By the previous considerations we conclude that
\[
\begin{split}
  q(S) 
  &= \sum_{\chi\in\Ghat} h^1(X,K_X+L_\chi) \\ 
  &= \sum_{W\in\FFF}\sum_{U\in\Uu(W)}|U|_\nn^M
  \cdot h^1\big(\PP^2,\Oo_{\PP^2}(-3+h_\cC(W))\otimes\Jj(U\cdot \cC)\big),
\end{split}
\]
where 
\[
  |U|_\nn^M
  = \card\Bigg\{\!
  (a^1,\ldots,a^s) \,\Big|\,
  0\leq a^j<n_j, 
  \Bigg(\!\!\fpart{\sum_ja^j\mu_j^1/n_j},\ldots,
  \fpart{\sum_ja^j\mu_j^t/n_j}\!\!\Bigg)
  \in U
  \Bigg\}.
\]

\paragraph{Second case} If $H_\infty$ is not transverse to $\sum_iC_i$
and $S'$ is ramified above $H_\infty$, then,
supposing for simplicity that there is only one singular point $Q$ of
$\sum_iC_i$ lying on $H_\infty$, we have
\[
  \mu^\ast (\sum_iC_i+H_\infty) 
  = \sum_i\Ctilde_i+\Htilde_\infty
  +\sum_{P\neq Q}\ee^P\cdot\eE_P
  +(\ee_Q+\ee_\infty^Q)\cdot\eE_Q,
\]
with $\ee^P=\sum_i\ee_i^P$ and $\ee_\infty^Q=(1,\ldots)$.  
As in the first case, 
\begin{multline*}
  n_jL''_{\chi_j}
  \sim  \sum_i\mu_j^i\Ctilde_i + 
  \sum_{i,P}\mu_j^i\ee_i^P\cdot\eE_{P} \\
  + \Bigg(\!\ceil{\sum_i\mu_j^id_i/n_j}n_j-\sum_i\mu_j^id_i\Bigg)\Htilde_\infty
  + \Bigg(\!\ceil{\sum_i\mu_j^id_i/n_j}n_j-\sum_i\mu_j^id_i\Bigg)
  \ee_\infty^Q\cdot\eE_Q,
\end{multline*} 
hence 
\begin{multline*}
  L''_\chi \sim
  \sum_j a^jL''_{\chi_j}
  -\sum_i\floor{\sum_j\frac{a^j\mu_j^i}{n_j}}\Ctilde_i
  -\floor{\sum_{j=1}^s\frac{a^j}{n_j}
  \Big(\!\ceil{\sum_i\mu_j^id_i/n_j}n_j-\sum_i\mu_j^id_i\Big)}
    \Htilde_\infty \\
  -\sum_{P\neq Q}\floor{\sum_{i,j}\frac{a^j\mu_j^i}{n_j}\ee_i^P}\cdot\eE_P 
  -\floor{\sum_{i,j}\frac{a^j\mu_j^i}{n_j}\ee_i^Q+
  \sum_{j=1}^s\frac{a^j}{n_j}
  \Big(\!\ceil{\sum_i\mu_j^id_i/n_j}n_j-\sum_i\mu_j^id_i\Big)
  \ee_\infty^Q}\cdot\eE_Q.
\end{multline*}
By (\ref{eq:pull-back}) and 
$\Htilde_\infty\sim\Htilde-\ee_\infty^Q\cdot\eE_Q$,
\[
\begin{split}
  L_\chi 
  &\sim 
  \Bigg(\!
  \ceil{\sum_{i,j}\frac{a^j\mu_j^id_i}{n_j}}
  -\sum_i\floor{\sum_j\frac{a^j\mu_j^i}{n_j}}d_i
  \Bigg)\Htilde \\
  &\quad - \sum_{P\neq Q}
  \Bigg(\!\floor{\sum_{i,j}\frac{a^j\mu_j^i}{n_j}\ee_i^P}
  -\sum_{P,i}\floor{\sum_j\frac{a^j\mu_j^i}{n_j}}\ee_i^P
  \Bigg)\cdot\eE_P \\
  &\quad - \Bigg(\!
  \floor{\sum_{i,j}\frac{a^j\mu_j^i}{n_j}\ee_i^Q
    -\sum_{i,j}\frac{a^j\mu_j^id_i}{n_j}\ee_\infty^Q}
  -\sum_i\floor{\sum_j\frac{a^j\mu_j^i}{n_j}}\ee_i^Q
  +\ceil{\sum_{i,j}\frac{a^j\mu_j^id_i}{n_j}}\ee_\infty^Q
  \!\Bigg)\cdot \eE_Q.
\end{split}
\]
Set $x^i=\fpart{\sum_ja^j\mu_j^i/n_j}$, $1\leq i\leq t$, as in
(\ref{eq:theX}).  Then
\begin{multline}
  \label{eq:nonTransverse}
  L_\chi \sim \ceil{\sum_{i=1}^td_ix^i}\Htilde 
  -\sum_{P\neq Q}\floor{\sum_{i=1}^tx^i\ee_i^P}\cdot\eE_P \\
  -\Bigg(\!\floor{\sum_ix^i\ee_i^Q-\sum_id_ix^i\ee_\infty^Q}
    +\ceil{\sum_id_ix^i}\ee_\infty^Q
  \!\Bigg)\cdot\eE_Q.
\end{multline}
In the formula for the irregularity, two things may happen.  Either
$\sum_id_ix^i $ is an integer and the superabundances involved can be
dealt with as before, or $\sum_id_ix^i $ is not an integer.  In this
latter situation set $C_0=H_\infty$, $d_0=1$, $\ee_0^P=\zero$ if
$P\neq Q$ and
\[
  x^0 = \ceil{\sum_id_ix^i}-\sum_id^ix_i.
\]
Then
\[
\begin{split}
  L_\chi 
  &\sim \bigg(\! x^0+\sum_{i=1}^td_ix^i \!\bigg)\Htilde 
  -\sum_{P\neq Q}\floor{\sum_{i=1}^tx^i\ee_i^P}\cdot\eE_P
  -\floor{x^0\ee_0^Q+\sum_{i=1}^tx^i\ee_i^{Q}}\cdot\eE_Q \\
  &= \bigg(\! \sum_{i=0}^td_ix^i \!\bigg)\Htilde 
  -\sum_P\floor{\sum_{i=0}^tx^i\ee_i^P}\cdot\eE_P
\end{split}
\]
and the formula for the irregularity follows as before 
replacing $\sum_{i=1}^tC_i$ by $\sum_{i=0}^tC_i$, $\cC$ by 
$\overline{\cC}=(C_0,\ldots,C_t)$ and $M$ by $\Mbar$, where 
\[
  \Mbar = 
  \begin{bmatrix}
    \mu_1^0 & \cdots & \mu_s^0 \\
    & M 
  \end{bmatrix}
\]
with
$\mu_j^0=\ceil{\sum_{i=1}^t\mu_j^id_i/n_j}n_j-\sum_{i=1}^t\mu_j^id_i$.
To end the proof it remains to show that
\[
  x^0 = \fpart{\sum_j\frac{a^j\mu_j^0}{n_j}}.
\]
But this is clear, since
\[
\begin{split}
  \fpart{\sum_j\frac{a^j\mu_j^0}{n_j}}
  &= \fpart{\sum_{j=1}^s\bigg(
    \ceil{\frac{1}{n_j}\sum_{i=1}^t\mu_j^id_i}n_j-\sum_{i=1}^t\mu_j^id_i
    \bigg)\frac{a^j}{n_j}} \\
  &= \fpart{\ceil{\sum_{i,j}d_i\frac{a^j\mu_j^i}{n_j}}
    -\sum_{i,j}d_i\frac{a^j\mu_j^i}{n_j}} \\
  &= \ceil{\sum_{i=1}^t d_i\fpart{\sum_{j=1}^s\frac{a^j\mu_j^i}{n_j}}}
    -\sum_{i=1}^t d_i\fpart{\sum_{j=1}^s\frac{a^j\mu_j^i}{n_j}}
\end{split}
\]
and this equals $x^0$ by the definition of the $x^i$ when 
$1\leq i\leq t$.
\qed

\section{Applications and examples}
\label{s:applications}

\subsection{Asymptotic behaviour of the irregularity}
In setting out to look for applications of Theorem~\ref{th:q} it seems
best to start with the asymptotic behaviour of the irregularity of the
abelian coverings of the projective plane described by E.~Hironaka in
\cite{Hiro}.

\begin{cor}
  \label{c:q}
Let $S'\to\PP^2$ be the $(\ZZ/n\ZZ)^s$-abelian covering defined 
by $L_{\chi_j}\sim C_j+(\ceil{d_j/n}n-d_j)H_\infty$, $1\leq j\leq s$,
with $d_j=\deg C_j$.  If $S$ is the normalization of $S'$, then 
\[
  q(S) = \sum_{W\in\FFF(\cC)}\sum_{U\in\Uu(W)}
  |U|_n \cdot h^1\big(\PP^2,\Oo_{\PP^2}(-3+h_{\cC}(W))
  \otimes\Jj(U\cdot\cC)\big)
\]
where the partition $\cC$ is given either by the curves
$C_1,\ldots,C_s$ or by the curves $H_\infty$, $C_1,\ldots,C_s$,
depending on whether or not $H_\infty$ is transverse to $\sum_iC_i$.
\end{cor}

\proof
By Theorem \ref{th:q}, when $\sum_iC_i$ is not transverse to
$H_\infty$ and the covering is branched along $H_\infty$,
\begin{equation} \label{eq:qbar}
  q(S) = 
  \sum_{W\in\FFF(\overline{\cC})}\sum_{U\in\Uu(W)}|U|_n^\Mbar
  \cdot h^1\big(\PP^2,\Oo_{\PP^2}(-3+h_{\overline{\cC}}(W))
  \otimes\Jj(U\cdot\overline{\cC})\big),
\end{equation}
where $\overline{\cC}=(H_\infty,\cC)$ and 
$\Mbar=\begin{bmatrix} \mu_1^0 & \cdots & \mu_s^0 \\& I_s \end{bmatrix}$,
with $\mu_j^0=\ceil{d_j/n}n-d_j$ for $1\leq j\leq s$.  Moreover,
$|U|_n^{\Mbar} = \card\phi\inv(U\cap[0,1)^{s+1})$, where the map 
$\phi : [0,1)^s\cap(1/n\ZZ)^s\to[0,1)^{s+1}$ associated to $\Mbar$ is
defined by
\[
  \phi\bigg(\frac{a^1}{n_1},\ldots,\frac{a^s}{n_s}\bigg) 
  = \bigg(\!\fpart{\sum_j\mu_j^0\frac{a^j}{n_j}},\frac{a^1}{n},\ldots,
    \frac{a^s}{n}\!\bigg).
\]
Since $\phi$ is injective, it follows that 
\[
  |U|_n^{\Mbar} = \card(U\cap[0,1)^{s+1}\cap(1/n\ZZ)^{s+1}) = |U|_n.
\]
\qed

Using Corollary~\ref{c:q} we can recover E.~Hironaka's result
concerning the asymptotic behaviour of the irregularity of the abelian
covering $\Sigma_n$.  See also \cite[Theorem~1.7]{Bu}, where N.~Budur
establish the quasi-polynomial behaviour of the Hodge numbers $h^{0,q}$
of the finite abelian coverings of a smooth $n$-dimensional variety.

\begin{cor}
Let $C\subset\PP^2$ be a reduced curve and
$\Omega=\PP^2\smallsetminus(C\cup H_\infty)$.  Let $\Sigma_n$ be the
unique normal abelian covering associated to the natural epimorphism
\[
  \pi_1(\Omega) \lra H_1(\Omega,\ZZ) \lra H_1(\Omega,\ZZ/n\ZZ) 
  \simeq (\ZZ/n\ZZ)^s,
\]
with $s$ the number of connected components of $C$.  Then
$q(\Sigma_n)$ is a quasi-polynomial function of $n$ of degree $\leq s$.
\end{cor}

\begin{definition*}
A function $f:\NN\to\NN$ is called a quasi-polynomial function if
there exists an integer $N>0$ and polynomials $P_0,\ldots,P_{N-1}$
such that $f(n)=P_j(N)$ if $n\equiv j\mod N$.
\end{definition*}

\proof
By Remark~\ref{r:hironaka}, the surface $\Sigma_n$ coincides with the
normalization of the abelian covering of the projective plane with
group $(\ZZ/n\ZZ)^s$, associated to $C=\sum_{j=1}^sC_j$ and $H_\infty$
and determined by
\[
  nL_j \sim C_j + (\ceil{d_j/n}n-d_j)H_\infty.
\]
By Corollary~\ref{c:q} we have
\[
  q(\Sigma_n) = 
  \sum_{W\in\FFF(\cC)}\sum_{U\in\Uu(W)}|U|_n
  \cdot h^1\big(\PP^2,\Oo_{\PP^2}(-3+h_{\cC}(W))
  \otimes\Jj(U\cdot\cC)\big),
\]
with $\cC$ the partition of $C$ induced by the curves $C_j$.  The
closure of the subset $U\subset W$ in $W$ represents a convex polytope
and its border in $W$ a finite union of convex polytopes.  Now, if
$\Pp\subset\RR^s$ is a convex polytope, the {\it Ehrhart
quasi-polynomial} of $\Pp$ is the function defined by
\[
  i(\Pp,n) = \card(n\Pp\cap\ZZ^s),
\]
where $n\Pp=\{n\xx\mid \xx\in\Pp\}$.  Clearly, the number $i(\Pp,n)$
is equal to the number of rational points in $\Pp\cap (1/n\ZZ)^s$.  We
refer the reader to \cite[Theorem 4.6.25]{St} where it is shown that
$i(\Pp,n)$ is indeed a quasi-polynomial whose degree is $\dim\Pp$.
The result follows.
\qed

\subsection{Cyclic coverings}

As a particular case of Theorem~\ref{th:q} we obtain the formula for
the irregularity of cyclic multiple planes.  This study has been
initiated by O.~Zariski in \cite{Za2} where he computed the
irregularity in case the branching curve has only nodes and cusps as
singularities.  Various generalizations have since been proposed to
Zariski's formula in \cite{Es,Li1,Li2,Li6,Va,Na}.

\begin{cor}
Let $C\subset\PP^2$ be a curve of degree $d$ and $H_\infty$ a line
transverse to $C$.  If $S_n$ is the normalization of the standard
cyclic $\ZZ/n\ZZ$-covering of the plane defined by the linear
equivalence $nL_\chi \sim C+(\ceil{d/n}n-d)H_\infty$, then
\[
  q(S_n) = \sum_{\substack{\xi\text{ jumping number of }C\\
      \xi\in1/(n\wedge d)\,\ZZ,\,\,0<\xi<1}}
  h^1(\PP^2,\Oo_{\PP^2}(-3+\xi d)\otimes\Jj(\xi\cdot C)).
\]
\end{cor}

\proof
Since $H_\infty$ is transverse to $C$, then the faces in the formula
(\ref{eq:q}) are points $\xi$ in the open interval $(0,1)$
corresponding to the jumping numbers of $C$ such that 
$\xi\cdot\deg C\in\ZZ$.  Moreover $|\xi|_n = \card(\{\xi\}\cap1/n\ZZ)$
equals $1$ or $0$ depending on weather $\xi n\in\ZZ$.  The result
follows.
\qed

In the case of a cyclic covering, if $H_\infty$ is not transverse to
$C$, then the faces in the formula (\ref{eq:qbar}) live in $\RR^2$
with euclidean coordinates $x$ and $x^\infty$.  They are of two types:
1) those for which $x^\infty=0$, in which case they are jumping
numbers for $C$ and the corresponding term in (\ref{eq:q}) is
determined as in the above corollary; 2) those for which
$x^\infty\neq0$ and the face is determined by an equation whose
homogeneous part must coincide with the linear form $h_{(H_\infty,C)}$
modulo the multiplication by a non-zero rational.  It may be said that
the results obtained for the irregularity are qualitatively different.
In the transverse situation the irregularity is constant as a function
of $n$.  In the non transverse situation the irregularity depends on
$n$.  More precisely, using Hironaka's result, it is a
quasi-polynomial of degree $\leq 1$.  The next example illustrates
this behaviour of the irregularity in the non transverse situation.

\begin{exa}
  \label{ex:cyclic}
The two ellipses $\Gamma_i$, $i=1,2$, with common tangents at $P$ and
$Q$ considered in Example~\ref{exa:2tacNodes} provide cyclic coverings
with maximal degrees for the quasi-polynomials that represent the
irregularity, whatever the relative position of the line at infinity.

If $H_\infty$ is transverse to $B=\Gamma_1+\Gamma_2$, then the 
$\ZZ/n\ZZ$-cyclic covering $S_n$ has non vanishing irregularity if and 
only if $n$ is divisible by $4$, since $3/4$ is the only jumping 
number of $B$, in which case 
\[
  q(S_n) = h^1(\PP^2,\Ii_{P,Q}) = 1.
\]
Now, if $H_\infty$ is the line through $P$ and $Q$, then we have seen
in Example~\ref{exa:2tacNodes} that there are two jumping walls $W_3$
and $W_4$, and three distinguished faces, the previous two and the
point $W_0$, the intersection of $W_3$ with the coordinate plane
$x^\infty=0$.  Set $C=\Gamma_1+\Gamma_2+H_\infty$.  By
Corollary~\ref{c:q}, we get, for $n\geq 5$%
\footnote{For $n=4$ the covering is branched only along the two conics
  and the irregularity equals $1$.  For $n=3$, the formula for the non
  transverse intersection applies and the irregularity equals $2$.
  For $n=2$ the covering is branched again only along the conics and
  $q=0$.}, 
\[
\begin{split}
  q(S_n) &= |W_0|_n\, h^1(\PP^2,\Jj(W_0\cdot C)) +
  |W_3\smallsetminus W_0|_n\, h^1(\PP^2,\Jj(W_3\smallsetminus W_0\cdot C))\\
  &\quad + |W_4|_n\, h^1(\PP^2,\Oo(1)\otimes\Jj(W_4\cdot C)) \\
  &= |W_0|_n\, h^1(\PP^2,\Ii_{P,Q}) 
  + |W_3\smallsetminus W_0|_n\, h^1(\PP^2,\Ii_{P,Q})
  + |W_4|_n\, h^1(\PP^2,\Ii_Z(1)) \\
  &= |W_3|_n\, h^1(\PP^2,\Ii_{P,Q})
  + |W_4|_n\, h^1(\PP^2,\Ii_Z(1)),
\end{split}
\]
where $Z$ is the subscheme supported at $P$ and $Q$ and determined by
the points and the directions of the tangents to the two conics at
$P$ and $Q$.  Since 
\[
  |W_l|_n = \card\{(x,x^\infty)\mid 4x+x^\infty=l,
  \,0\leq x<1,\,0\leq x^\infty<1, x,x^\infty\in1/n\ZZ\},
\]
$l=3,4$, it follows that 
\[
  q(S_n) = \floor{\frac{n+1}{4}} + \floor{\frac{n+3}{4}}.
\]
\end{exa}

\begin{exa}
If in the previous example we consider the abelian covering $\Sigma_n$ 
of $\PP^2$ with group $\ZZ/n\ZZ\times\ZZ/n\ZZ$ and branched along 
$C=\Gamma_1+\Gamma_2+H_\infty$ with the partition 
$\cC=(\Gamma_1,\Gamma_2,H_\infty)$,
then the formula for the irregularity is the
same but the faces are defined in $\RR^3$ by
\[
  W_l = \{(x^1,x^2,x^\infty)\mid 2x^1+2x^2+x^\infty=l\}, 
\]
$l=3,4$.  Then
\[
  |W_3|_n = \frac{1}{2}\floor{\frac{n}{2}}
  \bigg(n+\ceil{\frac{n}{2}}-3\bigg)
  \textq{and}
  |W_4|_n = \frac{1}{2}\floor{\frac{n-1}{2}}
  \bigg(\floor{\frac{n-1}{2}}-1\bigg),
\]
hence $q(\Sigma_n)=(n-1)(n-2)/2$.
\end{exa}

\subsection{Line arrangements with only triple points}

Next we want to point out that the formula (\ref{eq:q}) simplifies in
case the branching curve is a line arrangement $\Aa$ with only triple
points.

\begin{notation*}
Let $W$ be a face.  The subarrangement $\Aa_W$ will denote the minimal
subarrangement of $\Aa$ determined by the points that contribute to
$W$.  This subarrangement is unique since all points are triple points.
\end{notation*}

\begin{thm}
  \label{th:triplePoints}
Let $\Aa=\bigcup_{j=1}^m H_j$ be a line arrangement in the projective
plane and let $H_\infty$ be a line either of $\Aa$ or transverse to
$\Aa$.  Let $s=m-1$ in the former case and $s=m$ in the latter.  If
$S$ is the normalization of the standard abelian covering associated
to $\Aa$, the line $H_\infty$ and the group $G\simeq(\ZZ/n\ZZ)^s$,
then
\begin{equation}
\label{eq:triplePoints}
  q(S) = 
  \sum_{W\in\FFF(\Aa)}\sum_{U\in\Uu(W)}|U|
  \cdot
  h^1\big(\PP^2,\Oo_{\PP^2}(-3+2/3\deg\Aa_W)\otimes\Ii_{Z_U}\big).
\end{equation}
\end{thm}

\proof
For each singular point $P$ of the arrangement, the configuration of
exceptional divisors is reduced to only one divisor $E_P$, with
$e_P^\Aa=3$.  Moreover, $e_j^P=1$ or $0$ depending on weather the line
$H_j$ passes through $P$ or not.  Since $2/3$ is the only jumping number
smaller than $1$ for a triple point, the only relevant value is $2$.
It follows that for any $P$, the elementary wall $W_P$ is given by
\[
  W_P=\{(x^1,\ldots,x^s)\mid e_P^1x^1+\ldots+e_P^sx^s=2\}.
\]
Now, let $W$ be a bounded face in the formula for the irregularity.
There exists a unique minimal subarrangement $\Aa_W$ determined by the
points contributing to $W$.  

\paragraph{Claim} $h_\Aa(W)=2/3\deg(\Aa_W)$. 

\noindent 
Indeed, let $I\subset\{1,\ldots,s\}$ such that 
$\Aa_W=\bigcup_{i\in I}H_i$. Since $h_\Aa$ is constant along $W$, the
equation $\sum_{i\in I}x^i=h_\Aa(W)$ is a linear combination of the
equations defining $W$ in $[0,1)^{|I|}$.  Using this linear combination 
on the free term and also evaluated for $x^i=1$ for every $i\in I$, 
the result follows. 

\medskip

To end the proof of the theorem, it is sufficient, for any
$U\in\Uu(W)$, to consider the subscheme $Z_U$ of points $P$ that are
among the triple points of $A_W$ and for which 
$\floor{\sum_{i\in I}e_P^ix^i}=2$.  
\qed


\begin{exa}[The Ceva arrangement $A_1(6)$]
Let $\Aa$ be the Ceva arrangement of degree $6$ with three double
points and four triple points.  Let $S$ be the normalisation of the
abelian covering of $\PP^2$ branched along $\Aa$ with
$H_\infty\subset\Aa$ and group $(\ZZ/n\ZZ)^5$.  Then
\[
  q(S) = \frac{5(n-2)(n-1)}{2}.
\]
These surfaces are introduced by F.~Hirzebruch in \cite{Hi}.  If
$n=5$, the irregularity was computed by M.-N.~Ishida in \cite{Is}.  The
general case was dealt with by A.~Libgober in \cite{Li}.

The sub-arrangements that may have a non-zero contribution in the
formula (\ref{eq:triplePoints}) are either the pencil sub-arrangement
$\Aa_P$ of a triple point $P$, or the arrangement $\Aa$.  Now,
\[
  h^1\big(\PP^2,\Oo_{\PP^2}(-3+2/3\,\deg\Aa_P)\otimes
  \Jj(2/3\cdot\Aa_P)\big)
  = h^1(\PP^2,\Ii_P(-1)) = 1
\]
and, if $Z$ denotes the support of the triple points, 
\[
  h^1\big(\PP^2,\Oo_{\PP^2}(-3+2/3\,\deg\Aa)\otimes
  \Jj(2/3\cdot\Aa)\big)
  = h^1(\PP^2,\Ii_Z(1)) = 1.
\]
So,
\[
\begin{aligned}
  q(S) 
  &= \sum_P|W(\Aa_P)|\cdot h^1(\PP^2,\Ii_P(-1))
  + |W(\Aa)|\cdot h^1(\PP^2,\Ii_Z(1))\\
  &= \sum_P|W(\Aa_P)|+|W(\Aa)|.
\end{aligned}
\]
$|W(\Aa_P)|$ counts in how many ways $2n$ can be written as a sum of
three integers that vary in $\{0,1,\ldots,n-1\}$.  Let us denote this
integer by $\sigma_3(2n)$.  It follows that
$|W(\Aa_P)|=\sigma_3(2n)=(n-2)(n-1)/2$.  As for $|W(\Aa)|$, it counts
the number of solutions of
\[
  \frac{1}{n}\sum_{j=1}^6 a^j = 4
  \textq{and} 
  \frac{1}{n}\sum_{H_j\ni P}a^j = 2
  \textql{for every $P$}.
\]
This means that $a^1+a^2+a^3=2n$ and that $a^i=a^j$ if and only if the
lines $H_i$ and $H_j$ intersect in a double point of $\Aa$.  Hence
$|W(\Aa)|=\sigma_3(2n)$.  The result follows.
\end{exa}

\subsection{The arrangement dual to the arrangement defined by the
  inflexion points of a smooth cubic}

Another arrangement considered in \cite{Hi} is the dual to the
arrangement defined by the inflexion points of a smooth cubic.  Let
$\Aa$ be such an arrangement. It has degree $9$ and twelve triple
points as only singularities.  In particular, each line contains four
triple points.

\begin{pro}
  \label{p:HirzII}
Let $S$ the normalization of the abelian covering of $\PP^2$ branched
along $\Aa$ with $H_\infty\subset\Aa$ and group $\ZZ/n\ZZ)^8$.  Then
$q(S)= 8(n-1)(n-2)-2\delta_{n\!\!\mod3}^0$, where $\delta_i^j$ denotes
the Kronecker symbol.
\end{pro}

\begin{rem*}
The case $n=5$ is treated in \cite{Is} and the general case in
\cite{Li}.  In this latter paper the formula for the irregularity is
$8(n-1)(n-2)$, lacking the corrective term in case $n$ is divisible by
$3$.   
\end{rem*}

\proof
By Theorem \ref{th:triplePoints} we have to study faces for which the
degree of the corresponding subarrangement $\Aa_W$ is divisible by
$3$, \ie equals $3$, $6$ or $9$.  In the first case, $\Aa_W$ is a
pencil subarrangement with a single triple point.  It will be denoted
$\Aa_P$, with $P$ the triple point.  If $\sigma_3(2n)$ is as before
the number of ways $2n$ can be written as a sum of three integers from
the set $\{0,1,\ldots,n-1\}$, then
\[
  \sum_P|W(\Aa_P)|\cdot h^1\big(\PP^2,\Ii_P(-1)\big)
  = 12 \sigma_3(2n).
\]

In the second case, $\Aa_W$ is a Ceva subarrangement and it is easy to
see that such a subarrangement cannot exist.  As for the last case,
there are different faces $W$ such that $\Aa_W=\Aa$.  Let $W$ be
determined by nine points among the twelve triple points---at least
nine points are needed so that the corresponding $h^1$ might be non
zero.  Since any ten among the twelve points impose independent
conditions on cubics as soon as the two remaining points lie on a line
of the arrangement, we infer that $W$ is determined by nine points
such that there is no line of the arrangement containing any two of
the remaining three points.  Hence, through each of these three points
pass three lines of the arrangement.  Now, if $Z$ is the union of the
nine points that determine $W$, then $h^1(\PP^2,\Ii_Z(3))=1$.  If
$H_1$, $H_2$ and $H_3$ are the three lines through one of the triple
points not in $Z$, then summing up the conditions for the points of
$Z$ lying on each of these three lines, we obtain that
\[
  6n =3a^1+\sum_{j=4}^9a^j  = 3a^2+\sum_{j=4}^9a^j
  = 3a^3+\sum_{j=4}^9a^j.
\] 
Hence $a^j$ is constant for the lines passing through each of the
three missing points.  Let $a(W)$, $a'(W)$ and $a''(W)$ be these three
constant values.  By the preceding equalities, 
\begin{equation}
  \label{eq:aaa}
  a(W)+a'(W)+a''(W) = 2n.
\end{equation}
It follows that 
\begin{multline*}
  q(S) = \sum_P|W(\Aa_P)|\cdot h^1\big(\PP^2,\Ii_P(-1)\big)
  + \sum_{\substack{W\text{ given}\\\text{by $9$ points}}}
  \sum_{U\in\Uu(W)}|U|\cdot h^1(\PP^2,\Jj(U\cdot \Aa)(3)) \\
  + \sum_{\substack{W\text{ given}\\\text{by $10$ points}}}
  \sum_{U\in\Uu(W)}|U|\cdot h^1(\PP^2,\Jj(U\cdot \Aa)(3))
  + |W(\Aa)|\cdot h^1(\PP^2,\Jj(2/3\cdot\Aa)(3)),
\end{multline*}
since if a face is defined by eleven points, using (\ref{eq:aaa}), it
will be defined by all twelve in fact.  Moreover, in the two middle
sums, $h^1(\PP^2,\Jj(U\cdot \Aa)(3))=h^1(\PP^2,\Ii_{Z(U)}(3))=1$.
Indeed, if $U\in\Uu(W)$ and $W$ is defined by a set $Z$ of nine triple
points as above, the subscheme $Z(U)$ is the union of $Z$ and of
either one or two more points, depending on the comparison of $3a(W)$,
$3a'(W)$ and $3a''(W)$ with $2n$.

In the hereafter lemma it is shown that
$h^1(\PP^2,\Jj(2/3\cdot\Aa)(3))=2$.  From the preceding
considerations and since there are exactly four groups of three
points such that there is no line of the arrangement containing any
two among the three points,
\begin{multline*}
  \sum_{\substack{W\text{ given}\\\text{by $9$ points}}}
  \sum_{U\in\Uu(W)}|U|\cdot h^1(\PP^2,\Ii_{Z(U)}(3)) 
  + \sum_{\substack{W\text{ given}\\\text{by $10$ points}}}
  \sum_{U\in\Uu(W)}|U|\cdot h^1(\PP^2,\Ii_{Z(U)}(3)) \\
  = \sum_{\substack{W\text{ given}\\\text{by $9$ points}}}
  \big(|W|-\delta_{n\!\!\!\!\mod3}^0|W(\Aa)|\big)
  = 4(\sigma_3(2n)-\delta_{n\!\!\!\!\mod3}^0).
\end{multline*}
The corrective term $\delta_{n\!\!\mod3}^0$ is given by the fact
that if $n$ is divisible by $3$, then the point in $W$ corresponding
to the case $a(W)=a'(W)=a''(W)=2n/3$ is to be considered in the face
$W(\Aa)$.  We conclude that  
\[
  q(S) 
  = 12\sigma_3(2n) + 4(\sigma_3(2n)-\delta_{n\!\!\!\!\mod3}^0)
  + 2\delta_{n\!\!\!\!\mod3}^0 
  = 8(n-1)(n-2)-2\delta_{n\!\!\!\!\mod3}^0.
\]
\qed

\begin{lem}
$h^1(\PP^2,\Jj(2/3\cdot\Aa))=2$. 
\end{lem}

\proof
Let $Z$ denotes the twelve triple points of $\Aa$.  We apply the
trace-residual exact sequence to the three lines $H_1$, $H_2$ and
$H_3$ that pass through one of the triple points.  Let $P$ and $P'$ be
the points not lying on these lines.  The exact sequences are
\[  
0 \lra \Ii_{\Res_{H_1}Z}(2) \lra 
  \Ii_Z(3) \lra \Oo_{\PP^1}(-1) \lra 0, 
\]
\[  
  0 \lra \Ii_{\Res_{H_2}(\Res_{H_1}Z)}(1)
  \lra \Ii_{\Res_{H_1}Z}(2)
  \lra \Oo_{\PP^1}(-1) \lra 0 
\]
and 
\[
  0 \lra \Ii_{P+P'} \lra 
  \Ii_{\Res_{H_2}(\Res_{H_1}Z)}(1)
  \lra \Oo_{\PP^1}(-2) \lra 0, 
\]
since $\deg\Res_{H_2}(\Res_{H_1}Z)=3$ and
$\Res_{H_3}\big(\Res_{H_2}(\Res_{H_1}Z)\big)=P\cup P'$.  It follows
that
\[
  h^1(\PP^2,\Ii_Z(3)) = h^1(\PP^2,\Ii_{\Res_{H_2}(\Res_{H_1}Z)}(1)) 
  = h^1(\PP^2,\Ii_{P\cup P'})+h^2(\PP^1,\Oo_{\PP^1}(-2)) = 2.
\]
\qed

\subsection{The arrangement associated to the Hesse pencil}
Let $\Aa$ be the line arrangement associated to the Hesse pencil.  It
is composed by the lines of the four singular fibres of the pencil
generated by a smooth elliptic curve and its Hessian.  It has degree
$12$, and twelve double points and nine quadruple points as
singularities.  The double points correspond to intersection points of
lines from the same singular fibre.

\begin{pro}
Let $S$ the normalization of the abelian covering of $\PP^2$ branched
along $\Aa$ with $H_\infty\subset\Aa$ and group $(\ZZ/n\ZZ)^{11}$.  Then
\[
  q(S) = \frac{(n-1)(61n^2+97n-378)}{6}.
\] 
\end{pro}

\proof
Theorem \ref{th:q} must be used here.  The relevant values for each
singular point of $\Aa$ are $2$ and $3$.  The distinguished faces
appearing in the formula for the irregularity are of the following
types:

1) $W_P$ associated to $(P,2)$ for each singular point $P$.  The
corresponding term in the right hand member of (\ref{eq:qbar}) equals
$\sigma_4(2n)$ since the conditions are
\[
  \frac{1}{n}\sum_{j=1}^{12} a^j = 2
  \quad\text{and}\quad 
  \frac{1}{n}\sum_{H_j\ni P}a^j = 2.
\]

2) $W_P$ associated to $(P,3)$ for each singular point $P$.  Here the
corresponding term equals $\sigma_4(3n)$.

3) $W_\Bb$, with $\Bb$ a Ceva subarrangement.  The face $W_\Bb$ is
associated to the singular points of $\Bb$ seen in $\Aa$, with relevant
value $2$ for each one of them.  There are $54$ such subarrangements,
one for each choice of two fibers and two by two components in each
fiber---such a choice determines the four triple points of $\Bb$.  As
for the terms corresponding to $W_\Bb$ in the formula for the
irregularity, let $H_1,\ldots,H_6$ be the lines of $\Bb$ and let
$H_7,\ldots,H_{10}$ be the remaining lines through its triple points.
Furthermore we suppose that $H_i$ and $H_j$ intersect in a double
point if and only if $i+j=7$.  The defining conditions of $W_\Bb$ are
the four equalities corresponding to the triple points:
\[
  \frac{1}{n}(a^1+a^2+a^3+a^7) = 2
  \quad\text{and so on, plus }\quad
  a^{11} = a^{12} = 0.
\]
The corresponding cohomology group in the formula (\ref{eq:q}) is non
trivial if and only if $h_\Aa(W_\Bb)=\sum_{j=1}^{10}a^j/n=4$.
Summing these four conditions for the four points, we get 
\[
  2\sum_{j=1}^6a^j + \sum_{k=7}^{10}a^k = 8n.
\]
We conclude that $W_\Bb$ must be defined by $a^7=\cdots=a^{12}=0$,
$a^1=a^6$, $a^2=a^5$ and $a^3=a^4$, and $a^1+a^2+a^3=2n$.  Hence
$|W_\Bb| = \sigma_3(2n)$.

4) $W$ defined by all nine singular points: six with relevant value $2$,
and the remaining three with relevant value $3$.  There are three
lines $H_j$ that do not pass through the points whose relevant value is
$3$ and do not intersect in a point.  There are $72$ such
possibilities, $\binom{4}{3}\cdot3\cdot3\cdot2$---$\binom{4}{3}$
choices for the fibres with distinguished components, $3$ choices for
the distinguished component of the first fibre, $3$ for the second and
$2$ for the third.  But, applying the trace-residual exact sequence
with respect to the three components, we see that $h^1$ vanishes.

5) $W$ defined by all nine singular points: six with relevant value $3$,
and the remaining three with relevant value $2$---the configuration
obtained from the preceding one by exchanging $2$ with $3$.  As
before, $h^1=0$ too. 

6) $W$ defined by all nine singular points with relevant value $2$.
Again $h^1$ does not vanish if and only if $h_\Aa(W)=6$.  Hence the
linear system defining $W$ becomes
\[
  \frac{1}{n}\sum_{j=1}^{12} a^j = 6
  \quad\text{and}\quad 
  \frac{1}{n}\sum_{H_j\ni P}a^j = 2
  \quad\text{for every singular point $P$.}
\]
Summing up the conditions imposed by the multiple points yields 
\[
  3\sum_{j=1}^{12} a^j = 9\cdot 2n,
\]
hence $\sum_{H_j\ni P}a^j=2n$ for every $P$.  But then
\[
  6n = \sum_{P\in H_{j_0}}\sum_{H_j\ni P}a^j
  = 3a^{j_0}+\sum_{\substack{H_j \text{ not a component of the fibre} \\
      \text{that contains } H_{j_0}}}a^j.
\]
Hence $a^j$ is constant along each special fibre of the Hesse pencil
and $|W|=\sigma_4(2n)$.

7) $W$ defined by all nine singular points with relevant value $3$.
Arguing as in the previous case, $h_\Aa(W)=9$ and hence we have
\[
  \frac{1}{n}\sum_{j=1}^{12} a^j = 9
  \quad\text{and}\quad 
  \frac{1}{n}\sum_{j\in\alpha}a^j = 3
  \quad\text{for every $\alpha$},
\]
and eventually $|W|=\sigma_4(3n)$.  For the computation of the $h^1$
in this case, the trace residual sequence gives
\[
  0 \lra 
  \Ii_{\sum P}(3) \stackrel{u}{\lra}
  \Ii_{\sum 2P}(6) \stackrel{r}{\lra}
  \Oo_E  \lra 0,
\]
where $E$ is a smooth cubic from the Hesse pencil.  Now
$h^0(\PP^2,\Ii_{\sum P}(3))=2$ and \linebreak
$h^0(\PP^2,\Ii_{\sum 2p_\alpha}(6))\geq3$, hence $h^0r$ is
surjective yielding that $h^1(\PP^2,\Ii_{\sum 2p_\alpha}(6))=2$.

Summing up,
\[
\begin{aligned}
  q(S) &= 9\cdot\sigma_4(2n)+9\cdot\sigma_4(3n)+54\cdot\sigma_3(2n)
  +\sigma_4(2n)+2\cdot\sigma_4(3n) \\
  &= 10\,\frac{(n-1)(5n^2-n-12)}{6}+11\,\frac{(n-1)(n-2)(n-3)}{6}+
  54\,\frac{(n-1)(n-2)}{2} \\
  &= \frac{(n-1)(61n^2+97n-378)}{6}.
\end{aligned}
\]
\qed

\subsection{General multiple planes}

The last example we would like to consider is one that makes use of
Theorem \ref{th:q} in its full generality.  Let $\Aa$ be the Ceva's
arrangement with the lines $C_1,\ldots,C_6$ such that $C_i$ and $C_j$
determine a double point if and only if $i+j=7$.  Let $S'$ be the
$(\ZZ/5\ZZ)^3$-abelian covering of $\PP^2$ defined by the reduced
building data
\[
\begin{aligned}
  5L_{\chi_1} &\sim 3C_2+C_3+C_6\\
  5L_{\chi_2} &\sim 2C_1+2C_2+C_4 \\
  5L_{\chi_3} &\sim C_1+3C_3+C_5.
\end{aligned}
\]
It is one of the examples considered by M.-N.~Ishida in \cite[\S
6]{Is}, with $q(S)=10$, where $S$ is the normalization of $S'$.  In
\cite{Is} it is shown that this surface is a quotient of the
Hirzebruch surface constructed as an $(\ZZ/5\ZZ)^5$-abelian covering of
the plane, by the group $(\ZZ/5\ZZ)^2$.  It also verifies $c_1^2=c_2$.
Moreover it is asserted that the surface is isomorphic to the one
constructed by M.~Inoue (see \cite{In}) from the elliptic
modular surface of level $5$.

Let us show how the irregularity might be computed using
Theorem~\ref{th:q}.  There are non-reduced components in the branch
locus and $C_\infty$ is taken to be $C_6$.  We have
\[
  q(S) = \sum_{P\text{ triple point}}|W(\Aa_P)|_5^M+|W(\Aa)|_5^M,
\]
where
\[
  M = 
  \begin{bmatrix}
    0 & 2 & 1\\
    3 & 2 & 0 \\
    1 & 0 & 3\\
    0 & 1 & 0\\
    0 & 0 & 1\\
    1 & 0 & 0
  \end{bmatrix}
\]
and $|W|_5^M=\card\phi\inv(W\cap[0,1)^6)$,
with $\phi:[0,1)^3\cap(1/5\ZZ)^3\to[0,1)^6$ defined by 
\[
  \phi\bigg(\frac{a^1}{5},\frac{a^2}{5},\frac{a^3}{5}\bigg)
  = \bigg(\!\fpart{\sum_{j=1}^3m_j^1\,\frac{a^j}{5}},\ldots,
    \fpart{\sum_{j=1}^3m_j^6\,\frac{a^j}{5}}\!\bigg).
\]
Here as before, $\Aa_P$ is the pencil subarrangement determined by the
triple point $P$.  An easy computation gives $|W(\Aa_P)|_5^M=2$ for
every triple point.  Furthermore, since the equations
\[
  \fpart{\frac{2a^2+a^3}{5}} = \fpart{\frac{a^1}{5}}, \quad
  \fpart{\frac{3a^1+2a^2}{5}} = \fpart{\frac{a^3}{5}}, \quad
  \fpart{\frac{a^1+3a^3}{5}} = \fpart{\frac{a^2}{5}}
\]
and 
$\fpart{a^1/5}+\fpart{a^2/5}+\fpart{a^3/5}=2$
lead to the only solutions $(2,4,4)$ and $(4,3,3)$, it follows that
$|W(\Aa)|_5^M=2$ also.  Hence the irregularity is $10$.

\appendix

\section{Technical result} 

In the proof of Theorem~\ref{th:q} we used two technical results that
enabled us to describe the reduced building data of the normalization
of a standard covering---which is also a standard covering (see
\cite[Corollary~3.1]{Pa})---in terms of the initial reduced building
data.  The first result was Proposition 3.2 in \cite{Na}.  The second
is somehow similar and deals with the fourth step in the normalization
algorithm presented in \cite{Pa}.  It is the step peculiar to the
abelian situation.

\begin{pro}
  \label{p:normalization4}
Let $X$ be smooth and let $\pi\colon Y\to X$ be a standard abelian
covering determined by the set of reduced building data $\Ll_{\chi_j}$
and $B_f$, $1\leq j\leq s$ and $f\in\FFF$.  Let $C$ be a multiplicity
$1$ component of both $B_f$ and $B_g$, \ie $B_f=C+R_f$ and
$B_g=C+R_g$.  After the normalization procedure has been applied to
$C$ and $Y'\to Y$ is the new surface, if
$\chi=\chi_1^{a_1}\cdots\chi_s^{a_s}$, then
\[
\begin{split}
  L'_\chi
  &\sim \sum_{j=1}^s a_jL_{\chi_j}
  -\floor{\sum_{j=1}^s a_j\left(\frac{f(\chi_j)^\bullet}{m_f}
      +\frac{g(\chi_j)^\bullet}{m_g}\right)}C \\
  &\hspace{1em}
  -\floor{\sum_{j=1}^s\frac{a_jf(\chi_j)^\bullet}{m_f}}R_f
  -\floor{\sum_{j=1}^s\frac{a_jg(\chi_j)^\bullet}{m_g}}R_g
  -\sum_{h\neq f,g}
  \floor{\sum_{j=1}^s\frac{a_jh(\chi_j)^\bullet}{m_h}}B_h.
\end{split}
\]
\end{pro}

\proof
Assume that $f:\Ghat\to\ZZ/m_f$ and that $g:\Ghat\to\ZZ/m_g$.  Let $d$
and $m$ be the greatest common divisor of, and respectively the
smallest common multiple of $m_f$ and $m_g$.  If
$\phi:\ZZ/m_f\times\ZZ/m_g\to\ZZ/m$ is defined by
$\phi(1,0)=m_g/d$ and $\phi(0,1)=m_f/d$, then set
$f':\Ghat\to\ZZ/m_{f'}$ the morphism defined by the composition
\[
\Ghat \stackrel{f\times g}{\lra} \ZZ/m_f\times\ZZ/m_g
\stackrel{\phi}{\lra} \Im\overline{\phi}
\stackrel{\iota}{\lra} \ZZ/m_{f'} 
\]
where $\overline{\phi}$ is the morphism $\phi\circ(f\times g)$ and
$\iota$ the isomorphism defined by $\iota(m/m_{f'})=1$.  The
normalization of $Y$ along $C$ is constructed by modifying the
covering data as follows:
\[
L'_\chi \sim 
\begin{cases}
  L_\chi -C, & \textqr{if} 
  \dfrac{f(\chi)^\bullet}{m_f}+\dfrac{g(\chi)^\bullet}{m_g} \geq 1\\
  L_\chi, & \text{otherwise}
\end{cases}
\]
and
\[
  B'_f\sim B_f-C,\quad
  B'_g\sim B_g-C,\quad
  B'_{f'}\sim B_{f'}+C,\quad
  B'_h \sim B_h \textq{for} h\neq f,g,f'.
\]
Applying these modifications to (\ref{eq:L_chi}) gives
\[
\begin{split}
  L'_\chi
  &\sim \sum_ja_jL'_{\chi_j}-\sum_{h\in\FFF}
  \floor{\sum_j\frac{a_jh(\chi_j)^\bullet}{m_h}}B_h \\
  &\sim {\sum_i}'a_i(L_{\chi_i}-C)+
  {\sum_k}''a_kL_{\chi_k}-
  \floor{\sum_j\frac{a_jf(\chi_j)^\bullet}{m_f}}R_f-
  \floor{\sum_j\frac{a_jg(\chi_j)^\bullet}{m_g}}R_g\\
  &\hspace{1em}
  -\floor{\sum_j\frac{a_jf'(\chi_j)^\bullet}{m_{f'}}}(B_{f'}+C)
  -\sum_{h\neq f,g,f'}
  \floor{\sum_j\frac{a_jh(\chi_j)^\bullet}{m_h}}B_h,
\end{split}
\]
where the sum $\sum'$ runs over those $i$'s for which
$f(\chi)^\bullet/m_f+g(\chi)^\bullet/m_g\geq1$ and $\sum''$ over the
other $k$'s.  To prove the result it is sufficient to show that 
\[
  {\sum_i}'a_i+\floor{\sum_j\frac{a_jf'(\chi_j)^\bullet}{m_{f'}}}
  = \floor{\sum_{j=1}^s a_j\left(\frac{f(\chi_j)^\bullet}{m_f}
        +\frac{g(\chi_j)^\bullet}{m_g}\right)}.
\]
But
\[
\begin{split}
  \floor{\sum_j\frac{a_jf'(\chi_j)^\bullet}{m_{f'}}}
  &= \floor{\sum_j\frac{a_j}{m_{f'}}
    \bigg(
    \frac{m_g}{d}f(\chi_j)^\bullet+\frac{m_f}{d}g(\chi_j)^\bullet
    \bigg)^\bullet\,\frac{m_{f'}}{m}} \\
  &= \floor{{\sum_i}'+{\sum_k}''}.
\end{split}
\]
Since for each $i$ in the first sum
\[
  \bigg(
    \frac{m_g}{d}f(\chi_i)^\bullet+\frac{m_f}{d}g(\chi_i)^\bullet
  \bigg)^\bullet =
  \frac{m_g}{d}f(\chi_i)^\bullet+\frac{m_f}{d}g(\chi_i)^\bullet-m
\]
and for each $k$ in the second
\[
  \bigg(
    \frac{m_g}{d}f(\chi_k)^\bullet+\frac{m_f}{d}g(\chi_k)^\bullet
  \bigg)^\bullet =
  \frac{m_g}{d}f(\chi_k)^\bullet+\frac{m_f}{d}g(\chi_k)^\bullet,
\]
the identity follows.
\qed

\bigskip

\begin{flushleft}
  Daniel {\sc Naie}\\
  D\'epartement de Math\'ematiques\\
  Universit\'e d'Angers\\
  F-40045 Angers\\
  France\\
  \small{Daniel.Naie@univ-angers.fr}
\end{flushleft}


\begin{thebibliography}{99}

\bibitem{Bu}
  {\sc N. Budur},
  {\it Unitary local systems, multiplier ideals, and polynomial
    periodicity of Hodge numbers} (arXiv:math/0610382).

\bibitem{Di}
  {\sc A. Dimca},
  {\it Singularities and topology of hypersurfaces}.
  Universitext. Springer-Verlag, New York, 1992.

\bibitem{Es}
  {\sc H.~Esnault},
  Fibre de Milnor d'un c\^one sur une courbe alg\'ebrique plane.
  {\it Invent.~Math. 68\/} (1982), 477--496.

\bibitem{FaJo}
  {\sc Ch. Favre, M. Jonsson}
  Valuations and multiplier ideals, 
  {\it J. Amer. Math. Soc. 18\/} (2005) no. 3, 655--684.

\bibitem{GrRe}
  {\sc H.~Grauert, R.~Remmert},
  Komplexe R\"aume.
  {\it Math. Ann. 136\/} (1958), 245--318.

\bibitem{Hiro}
  {\sc E.~Hironaka},
  Polynomial periodicity for Betti numbers of covering surfaces.
  {\it Invent. Math. 108\/} (1992), 289--321. 

\bibitem{Hi}
  {\sc F.~Hirzebruch},
  Arrangements of lines and algebraic surfaces.
  In {\it  Artin, M., Task,
    J. (eds.) Arithmetic Geometry}, vol. II. Boston, Birkh\"auser, 1983.

\bibitem{In}
  {\sc M. Inoue},
  Some new surfaces of general type. 
  {\it Tokyo J. Math.  17\/}  (1994),  no. 2, 295--319. 


\bibitem{Is}
  {\sc M.-N.~Ishida},
  The irregularities of Hirzebruch's examples of surfaces of general
  type with $c_1^2=3c_2$. 
  {\it Math. Ann. 262\/} (1983), 407--420. 

\bibitem{La}
  {\sc R.~Lazarsfeld},
  {\it Positivity in algebraic geometry}.
  A Series of Modern Surveys in Mathematics, Springer-Verlag, Berlin,
  2004. 

\bibitem{Ja}
  {\sc T.~J\"arviletho},
  {\it Jumping numbers of a simple complete ideal in a two-dimensional
    regular local ring}. Ph.~D.~Thesis, University of Helsinky (2007).

\bibitem{Li1}
  {\sc  A.~Libgober},
  Alexander polynomial of plane algebraic curves and cyclic multiple
  planes. 
  {\it Duke Math.~J. 49\/} (1982), 833--851.

\bibitem{Li2}
  {\sc A.~Libgober},
  Homotopy groups of the complements to singular hypersurfaces. 
  {\it Bull. Amer. Math. Soc. 13\/} (1985), 49--52.

\bibitem{Li}
  {\sc A.~Libgober},
 Characteristic varieties of algebraic curves.
 In {\it Applications of algebraic geometry to coding theory, physics
   and computation} (Eilat, 2001), 215--254, 
 NATO Sci. Ser. II Math. Phys. Chem., 36,
 Kluwer Acad. Publ., Dordrecht, 2001.  
(or AG/9801070)

\bibitem{Li5}
  {\sc A. Libgober},
  Hodge decomposition of Alexander invariants.
  {\it Manuscripta Math. \/107}  (2002), 251--269.

\bibitem{Li6}
  {\sc A.~Libgober},
  Lectures on topology of complements and fundamental groups,
  arYiv:math.AG/0510049 (2005).

\bibitem{Na}
  {\sc D. Naie},
  The irregularity of cyclic multiple planes after Zariski.
  {\it L'enseignement math\'ematique 53\/} (2007), 265-305.

\bibitem{Na2}
  {\sc D. Naie},
  Jumping numbers of a unibranch curve on a smooth surface.
  {\it Manuscripta mathematica \/} (2008), .

\bibitem{Pa}
  {\sc R. Pardini},
  Abelian covers of algebraic varieties.
  {\it J.~reine angew.~Math. 417\/} (1991), 191--213.

\bibitem{SmTh}
  {\sc K. E. Smith, H. M. Thompson},
  Irrelevant Exceptional Divisors for Curves on a Smooth Surface.
  {\it Algebra, geometry and their interactions},  245--254,
  Contemp. Math., 448, Amer. Math. Soc., Providence, RI, 2007. 

\bibitem{St}
  {\sc R. P. Stanley},
  {\it Enumerative combinatorics}, Vol. 1.
  Cambridge Univ. Press, Cambridge, 1997

\bibitem{Va}
  {\sc M. Vaqui\'e},
  Irr\'egularit\'e des rev\^etements cycliques des surfaces
  projectives non singuli\`eres.  
  {\it Amer.~J.~Math. 114\/} (1992), no. 6, 1187--1199.

\bibitem{Za2}
  {\sc O. Zariski},
  On the linear connection index of the algebraic surfaces $z^n=f(x,y)$.
  {\it Proceedings Nat.~Acad.~Sciences 15\/} (1929), 494-501.

\bibitem{Za1}
  {\sc O. Zariski},
  On the irregularity of cyclic multiple planes.
  {\it Ann.~of Math. 32\/} (1931), 485-511.
\end{thebibliography}
\end{document}